\documentclass[12pt]{article}
 \usepackage{psfig}
\usepackage{amssymb,amsmath,amstext,amsthm,amsfonts}
\usepackage[dvips]{graphicx}
 \usepackage[ansinew]{inputenc}

\title{Random perturbations of non-uniformly\\ expanding
maps\thanks{ Work partially supported by FCT through
{Centro de Matemática da Universidade do Porto.}}}

\author{José F. Alves, Vítor Araújo}

\date{\today}

\begin{document}

\newcommand{\mcup}{\mbox{$\bigcup$}}
\newcommand{\mcap}{\mbox{$\bigcap$}}

\def \RR {{\mathbb R}}
\def \ZZ {{\mathbb Z}}
\def \NN {{\mathbb N}}
\def \PP {{\mathbb P}}
\def \TT {{\mathbb T}}

\def \ra {\rightarrow }
 \def \wh {\widehat }
 \def \un{\underline }
 \def \dist {\mbox{dist}}
 \def \ov {\overline}
 \def \supp {\mbox{supp}\, }
 \def \wlim {\mbox{$w^*$-}\lim_{n\ra\infty}\, }
 \def \distp{\mbox{d}_\PP }

\def \al {\alpha } \def \be {\beta } \def \de {\delta }
\def \ga {\gamma } \def \ep {\epsilon } \def \vfi {\varphi
} \def \th {\theta } \def \si {\sigma }

 \def \cf {\mathcal{F}}
 \def \cm {\mathcal{M}}
 \def \cn {\mathcal{N}}
 \def \cq {\mathcal{Q}}
 \def \cp {\mathcal{P}}
 \def \cc {\mathcal{C}}
  \def \ch {\mathcal{H}}

\newcommand{\dem}{\begin{proof}}
\newcommand{\cqd}{\end{proof}}

\newcommand{\qand}{\quad\text{and}\quad}

\newtheorem{maintheorem}{Theorem}
\renewcommand{\themaintheorem}{\Alph{maintheorem}}
\newcommand{\cmt}{\begin{maintheorem}}
\newcommand{\fmt}{\end{maintheorem}}

\newtheorem{T}{Theorem}[section]
\newcommand{\ct}{\begin{T}}
\newcommand{\ft}{\end{T}}

\newtheorem{Corollary}[T]{Corollary}
\newcommand{\cco}{\begin{Corollary}}
\newcommand{\fco}{\end{Corollary}}

\newtheorem{Proposition}[T]{Proposition}
\newcommand{\cpr}{\begin{Proposition}}
\newcommand{\fpr}{\end{Proposition}}

\newtheorem{Lemma}[T]{Lemma}
\newcommand{\cle}{\begin{Lemma}}
\newcommand{\fle}{\end{Lemma}}

\newtheorem{Remark}[T]{Remark}
\newcommand{\cre}{\begin{Remark}}
\newcommand{\fre}{\end{Remark}}

\newtheorem{Definition}[T]{Definition}
\newcommand{\cd}{\begin{Definition}}
\newcommand{\fd}{\end{Definition}}

\maketitle

\begin{abstract} We give both sufficient conditions and
necessary conditions for the stochastic stability of
non-uniformly expanding maps either with or without
critical sets. We also show that  the number of probability
measures describing the statistical asymptotic behaviour of
random orbits is bounded by the number of SRB measures if
the noise level is small enough.

 As an application
of these results we prove the stochastic stability of certain
classes of non-uniformly expanding maps introduced in
\cite{V} and \cite{ABV}. \end{abstract}

\section{Introduction}

In broad terms, Dynamical Systems theory is mostly
interested in describing the typical behaviour of orbits as
time goes to infinity, and understanding how this behaviour
is modified under small perturbations of the system. This
work concentrates in the study of the latter problem from a
probabilistic point of view.

Given a map $f$ from a manifold $M$ into itself, let
$(x_n)_{n\ge1}$ be the orbit of a given point $x_0\in M$,
that is $x_{n+1}=f(x_n)$ for every $n\ge 1$. Consider the
sequence of time averages of Dirac measures $\delta_{x_j}$
along the orbit of $x_0$ from time $0$ to $n$. A special
interest lies  on the study of the convergence of such time
averages for a ``large" set of points $x_0\in M$ and the
properties of their limit measures. In this direction, we
refer the work of Sinai \cite{Si} for Anosov
diffeomorphisms, later extended by Ruelle and Bowen
\cite{BR,R} for Axiom A diffeomorphisms and flows. In the
context of systems with no uniform hyperbolic structure
Jakobson \cite{J} proved the existence of such measures for
certain quadratic transformations of the interval
exhibiting chaotic behaviour. Another important
contribution on this subject was given by Benedicks and
Young \cite{BY}, based on the previous work of Benedicks
and Carleson \cite{BC1,BC2}, where this kind of measures
were constructed for Hénon two dimensional maps exhibiting
strange attractors. The recent work of Alves, Bonatti and
Viana \cite{ABV} shows that such measures exist in great
generality for systems exhibiting some non-uniformly
expanding behaviour.

The notion of stability that most concerns us can be
formulated in the following way. Assume that, instead of
time averages of Dirac measures supported on the iterates
of $x_0\in M$, we consider time averages of Dirac measures
$\delta_{x_j}$, where at each iteration we take $x_{j+1}$
close to $f(x_j)$ with a controlled error. One is
interested in studying the existence of limit measures for
these time averages and their  relation to the analogous
ones for unperturbed orbits, that is, the stochastic
stability of the initial system.

Systems with some uniformly hyperbolic structure are quite
well understood and stability results have been established
in general, see  \cite{Ki1,Ki2} and \cite{Yo}. The
knowledge of the stochastic behaviour of systems that do
not exhibit such uniform expansion/contraction is still
very incomplete. Important results on this subject  were
obtained by Benedicks, Young \cite{BY}, Baladi and Viana
\cite{BaV} for certain quadratic maps of the interval.
Another important contribution is the announced work of
Benedicks and Viana for Hénon-like strange attractors. As
far as we know these are the only results of this type for
systems with no uniform expanding behaviour.

In this work we present both sufficient conditions and
necessary conditions for the stochastic stability of
non-uniformly expanding dynamical systems.
 As an application of these results we prove
that the classes of non-uniformly expanding maps introduced
in \cite{V} and \cite{ABV} are stochastically stable.

\subsection{Statement of results}
 Let $f:M\ra M$ be a smooth map defined on a compact riemannian
 manifold $M$.
 We fix some normalized
 riemannian volume form $m$ on $M$ that we call {\em Lebesgue measure}.

 Given $\mu$ an $f$-invariant Borel probability measure on
 $M$, we say that $\mu$ is an \emph{SRB measure}
 if, for a positive
Lebesgue measure set of points $x\in M$, the averaged
sequence of Dirac measures along the orbit $(f^n(x))_{n\geq
0}$ converges in the weak$^*$ topology to $\mu$, that is,
 \begin{equation}\label{average}
 \lim_{n\ra +\infty}
 \frac{1}{n}\sum_{j=0}^{n-1}\vfi\big(f^n(x)\big)
 =\int \vfi\, d\mu
\end{equation} for every continuous map $\vfi:M\to \RR$. We
define the {\em basin}  of $\mu$ as the set of those points
$x$ in $M$ for which (\ref{average}) holds for all
continuous $\vfi$. The maps to be considered in this work
will only have a finite number of SRB measures whose basins
cover the whole manifold $M$, up to a set of zero Lebesgue
measure.

We are interested in studying random perturbations of the
map $f$. For that, we take a continuous map
 \[
 \begin{array}{rccl}
 \Phi:& T &\longrightarrow&  C^2(M,M)\\
 & t &\longmapsto & f_t
 \end{array}
 \]
from a metric space $T$ into  the space of $C^2$ maps from
$M$ to $M$, with $f=f_{t^*}$ for some fixed $t^*\in T$.
Given $x\in M$ we call the sequence $\big(f_{\un
t}^n(x)\big)_{n\ge1}$ a \emph{random orbit} of $x$, where
$\un t$ denotes an element $(t_1,t_2,t_3,\ldots)$ in the
product space $T^{\NN}$ and
 $$
 f^n_{\un t}= f_{t_n}\circ \cdots \circ
 f_{t_1}\quad \mbox{for}\quad n \ge1.
 $$
We also take a family $(\th_\ep)_{\ep>0}$ of probability
measures on $T$ such that $(\supp\th_\ep)_{\ep>0}$ is a
nested family of connected compact sets and
$\supp\th_\ep\rightarrow \{t^*\}$ when $ \ep\to 0$.
 We will also assume some quite general nondegeneracy
 conditions on  $\Phi$ and $(\th_\ep)_{\ep>0}$ (see the beginning of Section
 \ref{s.stationary}) and  refer to $\{\Phi,(\th_\ep)_{\ep>0}\}$
 as a {\em random perturbation} of $f$.

In the context of random perturbations of a  map we say
that a Borel probability measure $\mu^\ep$ on $M$ is
\emph{physical } if for a positive Lebesgue measure  set of
points $x\in M$,  the averaged sequence of Dirac
probability measures $\delta_{f_{\un t}^n(x)}$ along random
orbits $\big(f_{\un t}^n(x) \big)_{n\geq 0}$ converges in
the weak$^*$ topology to $\mu^\ep$ for $\th_\ep^\NN$ almost
every $\un t\in T^\NN$. That is,
 \begin{equation}\label{pertaverage}
 \lim_{n\ra +\infty}
 \frac{1}{n}\sum_{j=0}^{n-1}\vfi\big(f_{\un t}^n(x)\big)
 =\int \vfi \,d\mu^\ep \quad\mbox{for all  continuous  $\vfi\colon
 M\rightarrow\RR$}
\end{equation}
 and $\th^\NN_\ep$ almost every $\un t\in
T^\NN.$
 We denote the set of points $x\in M$ for which
 (\ref{pertaverage}) holds by $B(\mu^\ep)$ and call it the
 \emph{basin of $\mu^\ep$}.
The map $f\colon M\ra M$ is said to be {\em stochastically
stable} if the weak$^*$ accumulation points (when $\ep>0$
goes to zero) of the physical probability measures of $f$
are convex linear combinations of the (finitely many) SRB
measures of~$f$.


\subsubsection{Local diffeomorphisms}

Let $f:M\ra M$ be a $C^2$
local diffeomorphism of the manifold $M$.
 We say that $f$ is {\em non-uniformly expanding} if there
is some constant $c>0$ for which
 \begin{equation} \label{liminf1}
\limsup_{n\ra +\infty}\frac{1}{n}
\sum_{j=0}^{n-1}\log\|Df(f^j(x))^{-1}\|\leq -c<0
 \end{equation}
 for
Lebesgue almost every $x\in M$. It was proved in \cite{ABV}
that for a non-uniformly expanding local diffeomorphism $f$
the following holds:
 \begin{itemize}
 \item[(P)]
 {\em There is a finite number of
ergodic absolutely continuous (\emph{SRB}) $f$-invariant
probability measures $\mu_1,\dots ,\mu_p$ whose basins
cover a full Lebesgue measure subset of $M$. Moreover,
every absolutely continuous $f$-invariant probability
measure $\mu$ may be written as a convex linear
combination of $\mu_1,\dots ,\mu_p$: there are real numbers
$w_1,\dots,w_p\geq 0$ with $w_1+\cdots+w_p=1$ for which
$\mu=w_1\mu_1+\cdots+w_p\mu_p$}.
 \end{itemize}
The proof of the previous result was based on the
existence of  $\al$-hyperbolic times
for the points in $M$:
given $0<\al<1$, we say that $n\in\ZZ^+$ is a {\em
$\al$-hyperbolic time} for the point $x\in M$ if
 \begin{equation}\label{htimesc}
 \prod_{j=n-k}^{n-1}
\|Df(f^j(x))^{-1}\|\leq \al^k \quad\mbox{for every}\quad
1\leq k\leq n.
 \end{equation}
The existence of (a positive frequency of) $\al$-hyperbolic
times for  points $x\in M$ is a consequence of the
hypothesis of non-uniform expansion of the map $f$ and
permits us to define a map $h:M\to\ZZ^+$ giving the first
hyperbolic time for $m$ almost every $x\in M$.

In the context of random perturbations of a non-uniformly
expanding
 map we are also able to prove a result on the finitness
of physical measures.

\cmt \label{t.finite1}
 Let $f\colon M\rightarrow M$ be a  $C^2$ non-uniformly expanding
 local diffeomorphism. If $\ep>0$ is sufficiently small, then there
 are physical measures  $\mu^\ep_1,\dots, \mu^\ep_l$
 (with $l$ not depending on
 $\ep$)
such that:
 \begin{enumerate} \item  for each  $x\in M$ and
$\th_\ep^\NN$ almost every $\un t\in T^\NN$, the average of
Dirac measures $\delta_{f_{\un t}^n(x)}$ converges in the
weak$^*$ topology to some  $\mu^\ep_i$ with $1\le i\le l$;
\item for each $1\le i\le l$ we have
 $$
 \mu^\ep_i=\wlim\frac1n\sum_{j=0}^{n-1}\int\big(f_{\un t}^j\big)_*
 \big(m\mid B(\mu^\ep_i)\big)\,d\th_\ep^\NN(\un t),
 $$
 where $m\mid B(\mu^\ep_i)$ is the normalization of the
 Lebesgue measure restricted to $B(\mu^\ep_i)$;
 \item   if $f$ is
 transitive, then $l=1$.
 \end{enumerate}
 \fmt

We say that the map $f$ is {\em non-uniformly expanding for
random orbits} if there is some constant $c>0$ such that
for $\ep>0$ small enough
 \begin{equation} \label{liminf2}
\limsup_{n\ra +\infty}\frac{1}{n}
\sum_{j=0}^{n-1}\log\|Df(f_{\un t}^j(x))^{-1}\|\leq -c<0,
 \end{equation}
for $\th_\ep^\NN\times m$ almost every $(\un t,x)\in
T^\NN\times M$. Similarly to the deterministic situation,
condition (\ref{liminf2}) permits us to introduce a notion
of $\al$-hyperbolic times for points in $T^\NN\times M$
 and define a map
 $$
 h_\ep\colon T^\NN\times M\ra \ZZ^+
 $$
by taking $h_\ep(\un t,x)$ the first $\al$-hyperbolic time
for the point $(\un t,x)\in T^\NN\times M$ (see
Section~\ref{s.distortion}). Assuming that
 $h_\ep$ is integrable with respect to $\th^\NN_\ep\times
 m$,
 then
 \begin{equation}\label{c.unif}
 \|h_\ep\|_1=\sum_{k=0}^{\infty}k\,
 (\th_\ep^\NN\times m)
 \big(\{(\un t,x)\colon h_\ep(\un t,x)=k\}\,\big)
 <\infty.
 \end{equation}
We say that the family  $(h_\ep)_{\ep>0}$ has {\em uniform
$L^1$-tail}, if the series in (\ref{c.unif}) converges uniformly
to $\|h_\ep\|_1$ (as a series of functions on the variable $\ep$).

\cmt \label{t.stc}
 Let $f\colon M\rightarrow M$ be a non-uniformly expanding $C^2$
  local diffeomorphism.
 \begin{enumerate}
 \item If $f$ is stochastically stable,
 then $f$ is non-uniformly expanding for random orbits.
 \item If $f$ is non-uniformly expanding for random orbits
 and $(h_\ep)_\ep$ has uniform $L^1$-tail,
then $f$ is stochastically stable.
 \end{enumerate}
\fmt


\subsubsection{Maps with critical sets}

Similar results to those presented for random perturbations
of local diffeomorphisms will also be obtained for maps
with  critical sets in the sense of~\cite{ABV}. We start by
describing the class of maps that we are going to consider. Let
$f\colon M\ra M$ be a continuous map of the compact manifold $M$
that fails to be a $C^2$ local diffeomorphism on a critical set $\cc\subset M$ with zero Lebesgue measure. We assume that $f$ {\em behaves like a power of the distance} close to
the critical set $\cc$:
there are constants $B>1$ and $\be>0$ for which

 \begin{itemize}
 \item[(S1)]
\hspace{.1cm}$\displaystyle{\frac{1}{B}\dist(x,\cc)^{\be}\leq
\frac{\|Df(x)v\|}{\|v\|}\leq B\dist(x,\cc)^{-\be}}$;
 \item[(S2)]
\hspace{.1cm}$\displaystyle{\left|\log\|Df(x)^{-1}\|-
\log\|Df(y)^{-1}\|\:\right|\leq
B\frac{\dist(x,y)}{\dist(x,\cc)^{\be}}}$;
 \item[(S3)]
\hspace{.1cm}$\displaystyle{\left|\log|\det Df(x)^{-1}|-
\log|\det Df(y)^{-1}|\:\right|\leq
B\frac{\dist(x,y)}{\dist(x,\cc)^{\be}}}$;
 \end{itemize}
 for every $x,y\in M\setminus \cc$ with
$\dist(x,y)<\dist(x,\cc)/2$ and $v\in T_x M$. Given
$\delta>0$ we define the $\delta$-{\em truncated distance}
from $x\in M$ to $\cc$
 $$ \dist_\delta(x,\cc)= \left\{
\begin{array}{ll} 1 & \mbox{if }\dist(x,\cc)\geq \delta,\\
\dist (x,\cc) & \mbox{otherwise.} \end{array} \right. $$

Assume that $f$ is a non-uniformly expanding map, in the
sense that there is $c>0$ such that the limit in
(\ref{liminf1}) holds for Lebesgue almost every $x\in M$
(recall that we are taking $\cc$ with zero Lebesgue
measure) and, moreover, suppose that the orbits of $f$ have
{\em slow approximation to the critical set}: given small
$\gamma >0$ there is $\delta >0$ such that
 \begin{equation}
 \label{limsup1}
\limsup_{n\ra +\infty}\frac{1}{n}
\sum_{j=0}^{n-1}-\log\dist_\delta(f^j(x), \cc)\leq \gamma
 \end{equation}
for Lebesgue almost every $x\in M$. The results in
\cite{ABV} show that in this situation we obtain the
same conclusion on the finiteness of SRB measures for such an $f$, also holding property (P).

In order to prove the stochastic stability of maps with
critical sets we need to restrict the class of
perturbations we are going to consider: we take maps $f_t$
with the same critical set $\cc$ and impose that
 \begin{equation}\label{e.perturbation}
 Df_t(x)=Df(x) \quad\mbox{for every $x\in M\setminus\cc$ and $t\in T$}.
 \end{equation}
This may be implemented, for instance, in parallelizable
manifolds (with an additive group structure, e.g. tori
$\TT^d$ or cylinders $\TT^{d-k}\times\RR^k$) by considering
$$T=\{t\in \RR^d\colon \|t\|\leq \ep_0\}$$
 for some
$\ep_0>0$, $\th_\ep$ the normalized Lebesgue measure on the
ball of radius $\ep\leq \ep_0$, and taking $f_t=f+t$; that
is, adding at each step a random noise to the unperturbed
dynamics.

For the case of maps with critical sets we also need to
impose an analog of condition~(\ref{limsup1}) for random
orbits; we assume {\em slow approximation of random orbits
to the critical set}: given any small $\gamma >0$ there is
$\delta >0$ such that \begin{equation}
 \label{limsup2}
\limsup_{n\ra +\infty}\frac{1}{n}
\sum_{j=0}^{n-1}-\log\dist_\delta(f^j_{\un t}(x), \cc)\leq
\gamma
 \end{equation}
for $\th_\ep^\NN\times m$ almost every $(\un t,x)\in
T^\NN\times M$ and  small $\ep>0$. Results  similar to
those presented for local diffeomorphisms on the finiteness
of physical measures  can also  be obtained in this case.

\cmt \label{t.finite2}
 Let $f\colon M\rightarrow M$ be
 a $C^2$ non-uniformly expanding  map behaving
 like a power of the distance close to the critical set
 $\cc$, and whose orbits have slow
 approximation to  $\cc$. If $f$ is non-uniformly expanding for
 random orbits and random orbits have slow
 approximation to  $\cc$,
 then we arrive at the
same conclusions of Theorem~\ref{t.finite1}.
 \fmt

The property of non-uniform expansion
 for random orbits, together with the slow approximation of random
orbits to the critical set permit us to introduce a notion
of $(\al,\delta)$-hyperbolic times for points in $(\un
t,x)\in T^\NN\times M$ and define a map
 $$
h_\ep\colon T^\NN\times M\ra \ZZ^+,
 $$
by taking $h_\ep(\un t, x)$ the first
$(\al,\delta)$-hyperbolic time for the point $(\un t,x)\in
T^\NN\times M$, see Section \ref{s.distortion}. Assuming
that $h_\ep$ is integrable with respect to $\th_\ep\times
m$, then we obtain an analog to (\ref{c.unif}), which
enables us to define a notion of {\em uniform $L^1$-tail}
exactly in the same way as before.

Due to the fact that $\log\|Df^{-1}\|$ is not a continuous map
(it is not even everywhere  defined) we are not able to
present in this context a similar to Theorem~\ref{t.stc} in
all its strength. However, we obtain the same kind of
conclusion of the second item of Theorem~\ref{t.stc}.

\cmt \label{t.sts}
 Let $f\colon M\rightarrow M$ be
 non-uniformly expanding $C^2$ map behaving like a power of the
 distance close to its critical set $\cc$ and whose orbits have slow
 approximation to $\cc$. Assume that $f$ is non-uniformly
expanding for random orbits and random orbits have slow
approximation to $\cc$. If $(h_\ep)_\ep$ has uniform
$L^1$-tail, then $f$ is stochastically stable. \fmt

As a major application of the previous theorem we are
thinking of a class of maps on the cylinder $S^1\times \RR$
introduced in~\cite{V}. Subsequent
works~\cite{Al} and~\cite{AV} showed that such systems are
topologically mixing (thus transitive) and have a unique
SRB measure. In Section \ref{s.applications} we prove that
Viana maps satisfy the hypotheses of Theorem~\ref{t.sts},
hence being stochastically stable. An application of
Theorem~\ref{t.stc} will also be given in
Section~\ref{s.applications} for an open class of local
diffeomorphisms introduced in~\cite[Appendix A]{ABV}.


\section{Distortion bounds} \label{s.distortion}

In this section we generalize some of the results in
\cite{Al} and \cite{ABV} for the setting of stochastic
perturbations of a non-uniformly expanding map. These
results  will be proved in the setting of maps with critical sets.
Then everything follows in the same way  for local
diffeomorphisms if we think of $\cc$ as being equal to
the empty set, with the only exception of a particular point
that we clarify in Remark \ref{r.hyp} below (due to
the fact that we are not assuming condition
(\ref{e.perturbation}) for  maps with no critical sets).
For the next definition we take
$0<b<\min\{1/2,1/(2\beta)\}$.

\cd Given $0<\al<1$ and $\delta>0$, we say that $n\in\ZZ^+$ is a
$(\al,\delta)$-hyperbolic time for $(\un t,x)\in T^\NN\times M$
if
 $$
 \prod_{j=n-k}^{n-1}
 \|Df_{t_{j+1}}(f^j_{\un t}(x))^{-1}\|\leq\al^k
 \quad\mbox{and}\quad
 \dist_\delta(f_{\un t}^{n-k}(x),\cc)\geq \al^{bk}
 $$
for every $1\leq k\leq n$.
\fd

The following lemma, due to Pliss \cite{Pli}, provides the
main tool in the proof of the existence of hyperbolic times
for points with non-uniform expansion on random orbits.

\cle \label{l.pliss} Let $H\ge c_2 > c_1 >0$ and
$\zeta={(c_2-c_1)}/{(H-c_1)}$. Given real numbers
$a_1,\ldots,a_N$ satisfying
 $$
 \sum_{j=1}^N a_j \ge c_2 N
\qand a_j\le H \;\;\mbox{for all}\;\; 1\le j\le N,
 $$
  there
are $l>\zeta N$ and $1<n_1<\ldots<n_l\le N$ such that $$
\sum_{j=n+1}^{n_i} a_j \ge c_1\cdot(n_i-n) \;\;\mbox{for
each}\;\; 0\le n < n_i, \; i=1,\ldots,l. $$
 \fle
  \dem See \cite[Lemma 3.1]{ABV}.
 \cqd

\cpr \label{p.hyp} There are $\al>0$ and $\delta>0$ for
which $\th_\ep^\NN\times m$ almost every $(\un t,x)\in
T^\NN\times M$ has some $(\al,\delta)$-hyperbolic time.
\fpr
 \dem
Let $(\un t,x)\in T^{\NN}\times M$ be a point
satisfying~(\ref{liminf2}). For large $N$ we have  $$
-\sum_{j=0}^{N-1} \log \left\| Df ( f^j_{\un t} (x) )^{-1}
\right\| \ge \frac{c}2 N >0, $$ by definition of
non-uniform expansion on random orbits. Fixing $\rho>\be$
we see that condition (S1) implies
\begin{equation}\label{ineq1} \left| \log \left\|
Df(x)^{-1}  \right\| \right| \le \rho \left| \log {\rm
dist\,}(x,\cc) \right| \end{equation} for every $x$ in a
neighborhood $V$ of $\cc$. Now we take $\ga_1>0$ so that
$\rho\ga_1\le c/10$ and let $\de_1>0$ be small enough to
get
 \begin{equation}\label{ineq2}
-\sum_{j=0}^{N-1} \log {\rm dist\,}_{\de_1} (f^j_{\un
t}(x), S) \le \ga_1 N\quad\mbox{for large $N$},
\end{equation} which is possible after
property~(\ref{limsup1}) of slow approximation to $\cc$.
 Moreover, fixing $H\ge \rho|\log \de|$
sufficiently large in order that it be also an upper bound
for for the set $\{ -\log \| Df^{-1}_{\un t}\| : t\in T, \;
x\in M\setminus V\}$, then the set
 $$
 E=\{ 1\le j\le N: -\log
\|Df(f^{j-1}_{\un t}(x))^{-1}\|>H \}
 $$
  is such that
$f^{j-1}_{\un t}(x)\in V$ for all $j\in E$ and
 $$
\rho\left|\log {\rm dist\,}(f^{j-1}_{\un t}(x),\cc)\right|
> -\log \left\| Df(f^{j-1}_{\un t}(x))^{-1} \right\| > H
\ge \rho|\log \de|
 $$
 i.e., ${\rm dist\,}(f^{j-1}_{\un
t}(x),\cc)<\de_1$, in particular
$\dist_{\de_1}(f^{j-1}_{\un t}(x),\cc)= \dist(f^{j-1}_{\un
t}(x),\cc)<\de_1$ for all $j\in E$. Hence, defining
 $$
a_j=\left\{ \begin{array}{lcr} -\log \left\|
Df(f^{j-1}_{\un t}(x))^{-1} \right\| & \mbox{if} & j\not\in
E \\ 0 & \mbox{if} & j\in E \end{array} \right. $$ it holds
$a_j\le H$ for $1\le j \le N$, and~(\ref{ineq1})
and~(\ref{ineq2}) imply $$ -\sum_{j\in E} \log \left\|
Df(f^{j-1}_{\un t}(x))^{-1}\right\| \le \rho\sum_{j\in E}
\left| \log \dist (f^{j-1}_{\un t}(x),\cc) \right| \le
\rho\ga_1 N. $$
 Since $\rho\ga_1\le c/10$ we deduce
  $$
\sum_{j=1}^N a_j = \sum_{j=1}^N \left( -\log \left\|
Df(f^{j-1}_{\un t}(x))^{-1} \right\| \right) - \sum_{j\in
E} \left( -\log \left\| Df(f^{j-1}_{\un t}(x))^{-1}
\right\| \right) \ge \frac25 c N.
 $$
By the previous arguments we may apply Lemma \ref{l.pliss}
to the sequence $a_j$ with $c_1=c/5$ and $c_2=2c/5$ (we may
suppose $H>c_1$ too by increasing $H$ if needed). Thus
there are $\zeta_1>0$ and $l_1>\zeta_1 N$ times $1\le q_1 <
\ldots < q_{l_1} \le N$ such that
\begin{equation}\label{ineq3} \sum_{j=n+1}^{q_i} -\log
\left\| Df(f^{j-1}_{\un t}(x))^{-1} \right\| \ge
\sum_{j=n+1}^{q_i} a_j \ge \frac{c}2(q_i-n) \end{equation}
for every $0\le n <q_i, \; i=1,\ldots,l_1$. We observe
that~(\ref{ineq3}) is just the first part of the
requirements on $(\al,\de)$-hyperbolic times for $(\un
t,x)$ if $\al=\exp(c/5)$.

Now we apply again Lemma~\ref{l.pliss}, this time to the
sequence $a_j=\log \dist_{\de_2} (f^{j-1}_{\un t}(x),\cc)$,
where $\de_2>0$ is small enough so that for $\ga_2>0$ with
$2\ga_2 (bc)^{-1}<\zeta_1$ we have by
assumption~(\ref{limsup1}) $$ \sum_{j=0}^{N-1} \log
\dist_{\de_2} (f^j_{\un t}(x) , \cc) \ge -\ga_2
N\quad\mbox{for large $N$}. $$ Defining $c_1=bc/2$, $c_2=
-\ga_2$, $H=0$ and $$ \zeta_2=\frac{c_2-c_1}{H-c_1}=
1-\frac{2\ga_2}{bc}, $$ Lemma~\ref{l.pliss} ensures that
there are $l_2\ge\zeta_2 N$ times $1\le r_1 < \ldots <
r_{l_2} \le N$ satisfying \begin{equation}\label{ineq4}
\sum_{j=n+1}^{r_i} \log \dist_{\de_2} (f^{j+1}_{\un t}(x),
\cc) \ge \frac{bc}2 (r_i -n ) \end{equation} for every
$0\le n < r_i$, $i=1,\ldots,l_2$. Let us note that the
condition on $\ga_2$ assures $\zeta_1+\zeta_2>1$. So if
$\zeta=\zeta_1+\zeta_2-1$, then there must be
$l=(l_1+l_2-N)\ge \zeta N$ and $1\le n_1 <\ldots <n_l \le
N$ for which~(\ref{ineq3}) and~(\ref{ineq4}) both hold.
This means that for $1\le i \le l$ and $1\le k \le n_i$ we
have
 $$
 \prod_{j=n_i-k}^{n_i} \left\| Df (f^j_{\un
t}(x))^{-1} \right\| \le \al^k \qand
\dist_{\de_2}(f^{n_i-k}_{\un t}(x),\cc)\ge \al^{bk},
 $$
 and
hence these $n_i$ are $(\al,\de)$-hyperbolic times for
$(\un t, x)$, with $\de=\de_2$ and $\al=\exp(c/5)$. It
follows that for $\th_\ep^\NN\times m$ almost every $(\un
t,x)\in T^\NN\times M$ there are (positive frequency of)
times $n\in\ZZ^+$ for which
 \begin{equation}\label{e.simple}
 \prod_{j=n-k}^{n-1} \|Df(f^j_{\un t}(x))^{-1}\|\leq\al^k
 \quad\mbox{and}\quad
 \dist_\delta(f_{\un t}^{n-k}(x),\cc)\geq \al^{bk}
  \end{equation}
for every $1\leq k\leq n$. Now the conclusion of the lemma
is a direct consequence of assumption (\ref{e.perturbation}).
\cqd

\cre\label{r.hyp}
In the setting of
random perturbations of a local diffeomorphism $f$
we may also derive
from the first part of (\ref{e.simple}) the existence
of hyperbolic times for $\th_\ep^\NN\times m$
almost every $(\un t,x)\in T^\NN\times M$ without assuming
condition (\ref{e.perturbation}).
Actually, let $(\un t,x)$ be a point in $ T^\NN\times M$ for which the
 first part of (\ref{e.simple}) holds. Taking the perturbations
$f_t$ in a sufficiently small $C^1$-neighborhood of $f$,  then
 $$
 \|Df_t(y)^{-1}\|\leq \frac{1}{\sqrt\al}\|Df(y)^{-1}\|
 $$
for every $y\in M$, which together with (\ref{e.simple}) gives
 $$
 \prod_{j=n-k}^{n-1} \|Df_t(f^j_{\un t}(x))^{-1}\|\leq
 \prod_{j=n-k}^{n-1} \frac{1}{\sqrt\al}\|Df(f^j_{\un t}(x))^{-1}\|
 \leq\al^{k/2}.
 $$
In the context of maps with no critical sets this $n$ may
be defined as a $\sqrt\al$-hyperbolic time for $(\un t,x)$
and all the results that we present below hold with  $\sqrt\al$-hyperbolic times replacing  $(\al,\delta)$-hyperbolic times for maps with critical sets.
 \fre

Proposition \ref{p.hyp} allows us to introduce a map
 $$
 h_\ep \colon T^\NN\times M\ra \ZZ^+,
 $$
by taking $h_\ep(\un t,x)$ as the first
$(\al,\delta)$-hyperbolic time for $(\un t,x)\in
T^\NN\times M$. We assume henceforth that the family
$(h_\ep)_{\ep>0}$ has uniform $L^1$-tail. For the next lemma we
fix $\de_1>0$ in such a way that
$4\delta_1<\min\{\de,\de^{\be}|\log \al|\}.$

\cle \label{l.contr}
 Given any $1\le j \le n$, we have
 \begin{equation*}\label{l.contr1}
  \|
Df(y)^{-1} \| \le \al^{-1/2} \| Df( f^{n-j}_{\un t}(x)
)^{-1} \|
 \end{equation*}
 for every  $y$ in the ball of radius $2\de_1\al^{j/2}$ around
$f^{n-j}_{\un t}(x)$.
 \fle
  \dem We are assuming
$\dist_{\de}(f^{n-j}_{\un t}(x),\cc)\ge \al^j$ since $n$ is
a $(\al,\de)$-hyperbolic time for $(\un t,x)$. This means
that
 $$
\dist (f^{n-j}_{\un t}(x),\cc)= \dist_{\de}(f^{n-j}_{\un
t}(x),\cc)\ge \al^{bj} \;\;\mbox{or else}\;\; \dist
(f^{n-j}_{\un t}(x),\cc)\ge \de.
 $$
  Either way it holds
$\dist(y,f^{n-j}_{\un t}(x)) \ge \dist (f^{n-j}_{\un
t}(x),\cc)/2$ because $b<1/2$ and $\de_1<\de/4<1/4$ for all
$y$ in the ball of radius $2\de_1\al^{j/2}$ around $
f^{n-j}_{\un t}(x)$. Therefore condition (S2) implies
 $$ \log\frac{\| Df(y)^{-1}
\|}{\|Df(f^{n-j}_{\un t}(x))^{-1}\|} \le
B\frac{\dist(f^{n-j}_{\un t}(x),y)} {\dist(f^{n-j}_{\un
t}(x),\cc)^\be} \le B\frac{2\de_1\al^{j/2}}{\min\{\al^{b\be
j}, \de^\be\}}.
 $$
  But $\al,\de <1$ and $b\be<1/2$ so
$\al^{j/2}<\al^{b\be j}$ and thus the right hand side of
the last expression is bounded from above by
$2B\de_1\de^{-\be}$. The assumptions on $\de_1$ assure this
last bound to be smaller than $\log \al^{-1/2}$, which
implies the statement. \cqd

\cpr\label{p.contr} There is $\de_1>0$ such that if $n$ is
$(\al,\delta)$-hyperbolic time for $(\un t,x)\in
T^\NN\times M$, then there is a neighborhood $V_n(\un t,x)$
of $x$ in $M$ such that \begin{enumerate} \item $f_{\un
t}^n$ maps $V_n(\un t,x)$ diffeomorphically onto the ball
of radius $\de_1$ around $f_{\un t}^n(x)$; \item for every
$1\leq k\leq n$ and $y,z\in V_k(\un t,x)$
 $$ \dist(f_{\un
t}^{n-k}(y),f_{\un t}^{n-k}(z)) \leq \al^{k/2}\dist(f_{\un
t}^{n}(y),f_{\un t}^{n}(z)).
 $$
  \end{enumerate}
   \fpr
   \dem
The proof will be by induction on $j\ge1$. First we show
that there is a well defined branch of $f^{-j}$ on a ball
of small enough radius around $f^j_{\un t}(x)$. Now we
observe that Lemma \ref{l.contr} gives for $j=1$
 $$
\| Df(y)^{-1} \| \le \al^{-1/2}\|Df(f^{n-1}_{\un
t}(x))^{-1}\| \le \al^{1/2},
 $$
 because $n$ is a
$(\al,\de)$-hyperbolic time for $(\un t,x)$. This means
that $f$ is a $\al^{-1/2}$-dilation in the ball of radius
$2\de_1\al^{1/2}$ around $ f^{n-1}_{\un t}(x)$.
Consequently there is some neighborhood $V_1(\un t, x)$ of
$f^{n-1}_{\un t}(x)$ inside the ball of radius
$2\de_1\al^{1/2}$ that is diffeomorphic to the ball of
radius $\de_1$ around $ f^n_{\un t}(x) $ through $f_{t_n}$,
when $f$ is a map with critical set
satisfying~(\ref{e.perturbation}).

For $j\ge1$ let us suppose that we have obtained a
neighborhood $V_j(\un t,x)$ of $f^{n-j}_{\un t}(x)$ such
that $f_{t_n}\circ\cdots\circ f_{t_{n-j+1}}\mid V_j(\un
t,x)$ is a diffeomorphism onto the ball of radius $\de_1$
around $ f^n_{\un t}(x) $ with
\begin{equation}\label{p.contr2}
\|Df(f_{t_{n-j+i+1}}\circ\cdots\circ
f_{t_{n-j+1}}(z))^{-1}\| \le
\al^{-1/2}\|Df(f^{n-j+i+1}_{\un t}(x))^{-1}\|
\end{equation} for all $z\in V_{j}(\un t,x)$ and $0\le
i<j$. Then, by Lemma \ref{l.contr} and under the assumption
that $n$ is a $(\al,\de)$-hyperbolic time for $x$,
 \begin{eqnarray*}
\big\| D\big(f_{t_n}\circ\cdots\circ
f_{t_{n-j}}(y)\big)^{-1} \big\| &\le& \prod_{i=0}^j \big\|
Df_{t_{n-j+i}}\big (f_{t_{n-j+i-1}}\circ\cdots\circ
f_{t_{n-j}}(y)\big)^{-1} \big\| \\ &\le& \prod_{i=0}^j
\al^{-1/2} \big\| Df_{t_{n-j+i}}\big(f^{n-j+i-1}_{\un
t}(x)\big)^{-1} \big\| \\ &\le&
(\al^{-1/2})^{j+1}\cdot\al^{j+1}=\al^{(j+1)/2}
 \end{eqnarray*}
  for every $y$ on the ball of radius
$2\de_1\al^{(j+1)/2}$ around $f^{n-j-1}_{\un t}(x)$ whose
image $f_{t_{n-j}}(y)$ is in $V_{j}(\un t,x)$ (above we
convention $f_{t_{n-j+i-1}}\circ\cdots\circ
f_{t_{n-j}}(y)=y$ for $i=0$).

This shows that  the derivative of $f_{t_n}\circ\cdots\circ
f_{t_{n-j}}$ is a $\al^{-(j+1)/2}$-dilation on the
intersection of
 $
f^{-1}_{t_{n-j}}\big(V_{j}(\un t,x)\big)
 $
 with the ball of radius  $2\de_1\al^{(j+1)/2}$ around
 $  f^{n-j-1}_{\un t}(x)$,
 and hence there is an
inverse branch of $f_{t_n}\circ\cdots\circ f_{t_{n-j}}$
defined on the ball of radius $\de_1$ around  $f^n_{\un
t}(x) $. Thus we may define $V_{j+1}(\un t,x)$ as the image
of the ball of radius $\de_1$ around $ f^n_{\un t}(x) $
under this inverse branch, and recover the induction
hypothesis for $j+1$. In this manner we get neighborhoods
$V_{j}(\un t, x)$ of $f^{n-j}_{\un t}(x)$ as above for all
$1\le j\le n$. \cqd

\cco\label{c.dist}
 There is a constant $C_1>0$ such that if
$\un t\in T^\NN$, $n$ is a $(\al,\de)$-hyperbolic time for
$x\in M$ and $y,z\in V_n(\un
t,x)$, then
 $$
 \frac{1}{C_1}\leq \frac{|\det Df_{\un
 t}^n(y)|}{|\det Df_{\un t}^n(z)|} \leq C_1.
 $$
 \fco
  \dem
For $1\le k\le n$ the distance between $f^k_{\un t}(x)$ and
either $f^k_{\un t}(y)$ or $f^k_{\un t}(z)$ is smaller than
$\al^{(n-k)/2}$ which is smaller than $\al^{b(n-k)}\le
\dist(f^k_{\un t}(x),\cc)$. So, by  (S3) we have
\begin{eqnarray*} \log\frac{|\det Df_{\un t}^n(y)|} {|\det
Df^n_{\un t}(z)|} &=& \sum_{k=0}^{n-1} \log \frac{|\det
Df_{t_{k+1}}(f^k_{\un t}(y))|} {|\det Df_{t_{k+1}}(f^k_{\un
t}(z))|}\\ & \le &\sum_{k=1}^{n-1}\log\frac{|\det
Df(f^k_{\un t}(y))|} {|\det Df(f^k_{\un t}(z))|} \\ &\le&
\sum_{k=0}^{n-1} 2B\frac{\al^{(n-k)/2}}{\al^{b\be(n-k)}},
\end{eqnarray*} and it is enough to take $C_1\le \exp
\left( \sum_{i=1}^\infty 2B \al^{(1/2-b\be)i} \right)$,
recalling that $b\be<1/2$ and also~(\ref{e.perturbation}).
\cqd


\section{Stationary measures} \label{s.stationary}

As mentioned before, we will assume the random
perturbations of the non-uniformly expanding map $f$
satisfy some {\em nondegeneracy conditions}: there
exists $0<\ep_0<1$ such that for every $0<\ep<\ep_0$ we may
take $n_0=n_0(\ep)\in\NN$ for which the following holds:

 \begin{enumerate}
 \item  there is $\xi=\xi(\ep)>0$ such
 that $ \left\{ f^n_{\un t}(x) \colon \un{t}\in
 (\supp\th_\ep)^{\NN} \right\} $ contains the ball of
 radius~$\xi$ around $f^n(x)$ for all $x\in M$ and
 $n\ge n_0$;
 \item
$(f^n_x)_*{\th}^{\NN}_\ep\ll m$  for all $x\in M$ and
 $n\ge n_0$.
 \end{enumerate}
 Here $(f^n_x)_*{\th}^{\NN}_\ep$
is the push-forward of ${\th}^{\NN}_{\ep}$ to $M$ via
$f^n_x:T^{\NN}\ra M$, defined as $f_x^n(\un t)=f_{\un t}^n(x)$.
 Condition~1 means that
perturbed iterates   cover a full neighborhood of the
unperturbed ones after a threshold for all sufficiently
small noise levels. Condition~2 means that sets of
perturbation vectors of positive ${\th}_{\ep}^{\NN}$
measure must send any point $x\in M$ onto subsets of $M$
with positive Lebesgue measure after a finite number of
iterates.

In~\cite[Examples 1 \& 2]{Ar1} it was shown that given any
smooth map $f:M\ra M$ of a compact manifold we can always
construct a random perturbation satisfying the
nondegeneracy conditions 1 and 2, if we take $T=\RR^p$,
$t^*=0$ and $\th_\ep$ is equal to the normalized
restriction of the Lebesgue measure to the ball of radius
$\ep$ around $0$, for a sufficiently big number $p\in\NN$
of parameters. For parallelizable manifolds the random
perturbations which consist in adding  at each step a
random noise to the unperturbed dynamics, as described in
the Introduction, clearly satisfy nondegeneracy conditions
1 and 2 for $n_0=1$.

\medskip

In the context of random perturbations of a map, we say that a set $A\subset M$ is {\em invariant} if $f_t(A)\subset A$, at least
for $t\in\supp(\th_\ep)$ with $\ep>0$ small.
The usual invariance of a measure with respect to a
transformation is replaced by the following one: a
probability measure $\mu$ is said to be \emph{stationary},
if for every continuous $\vfi:M\to\RR$ it holds
 \begin{equation}
 \label{eq.stationary}
 \int\vfi\,d\mu=
 \int\int\vfi\big(f_t(x)\big)\,d\mu(x)\,d\th_\ep(t).
 \end{equation}

\cre\label{re.accinvariant}
 If $(\mu^\ep)_{\ep>0}$ is a
family of stationary measures having $\mu_0$ as a weak$^*$
accumulation point when $\ep$ goes to $0$, then it follows
from~(\ref{eq.stationary}) and the convergence of
$\supp(\th_\ep)$ to $\{t^*\}$  that $\mu_0$ must be
invariant by $f=f_{t^*}$.
 \fre
It is not difficult to see (cf.~\cite{Ar1}) that a
stationary measure $\mu$ satisfies
 $$
 x\in\supp(\mu)
\quad\Rightarrow\quad f_t(x)\in\supp(\mu) \quad\mbox{for
all} \quad t\in\supp(\th_\ep)
 $$
 just by continuity of
$\Phi$.   This means that if $\mu$ is a
stationary measure, then $\supp(\mu)$ is an invariant set.
Nondegeneracy condition 1 ensures that the interior of
$\supp(\mu)$ is nonempty.

Let us write $\supp(\mu)$ as a disjoint union $\bigcup_i C_i$
of connected components and consider only those $C_i$ for
which $m(C_i)>0$ -- this collection is nonempty
since
$\supp(\mu)$ contains open sets. Moreover each $f_t$ must
permute these components for $t\in\supp(\th_\ep)$, because
$f_t(C_i)$ is connected by continuity,
$f_t(C_i)\subset\supp(\mu)$ by invariance, and
$m(f_t(C_i))>0$ since we have $(f_t)_*m\ll m$.

The connectedness of $C_i$ and continuity of $\Phi$ guarantee that
the abovemention perturbation of the components $C_i$
induced by $f_t$ does not depend on $t\in\supp(\th_\ep)$.
Indeed, supposing that $t,t^\prime\in\supp(\th_\ep)$ are
such that $$ f_t(C_i)\subset C_j \qand
f_{t^\prime}(C_i)\subset C_{j^\prime}, $$ then fixing some
$z\in C_i$ we have that $\{ f_t(z)\colon t\in\supp(\th_\ep)\}$
is a connected set intersecting both $C_j$ and
$C_{j^\prime}$ inside $\supp(\mu)$, and so $C_j =
C_{j^\prime}$.

We will show that these connected components are periodic
under the action induced by $f_t$ with $
t\in\supp(\th_\ep)$. After this, we may use nondegeneracy
condition 1 to conclude that each component contains a ball
of uniform radius and thus that each component satisfies
$m(C_i)> \mbox{const} >0$. Hence there existing only a
finite number of such components.

At this point it is useful to introduce the  skew-product map
\[
\begin{array}{rccc} F: & T^\NN\times M &\longrightarrow &
T^\NN\times M\\
 & (\un t, z) &\longmapsto & \big(\sigma(\un t),f_{t_1}(z)\big)
\end{array}
 \]
 where $\sigma$ is the left shift on
sequences $\un t=(t_1,t_2,\dots)\in T^\NN$.
 It is easy to check that the product measure $\th_\ep^\NN\times\mu$
is $F$-invariant, as so is the set
$\supp(\th_\ep^\NN\times\mu)=\supp(\th_\ep)^\NN\times\supp(\mu)$.

 \cle \label{pr.cycles}
 The support of a stationary measure
$\mu$ contains a finite number of connected components
arranged in cycles permuted by the action of $f_t$ for
$t\in\supp(\th_\ep)$.
 \fle
 \dem
Is is enough to obtain that each connected component $C_i$
is periodic under the action of $f_t$ for $t\in\supp(\th_\ep)$, in
the sense that $f^p_{\un t}(C_i)\subset C_i$ for some
$p\in\NN$ and all $\un t\in\supp(\th_\ep^\NN)$.
There are components $C_i$ with nonempty interior,
since the interior of $\supp(\mu)$ is nonempty.
So we may take a component
$C_i$ that contains some ball $B$.
Then
we have $m(B)>0$ and so
$(\th_\ep^\NN\times\mu)(\supp(\th_\ep^\NN)\times B)>0$.
Poincar\'e Recurrence Theorem now guarantees there is
$(\un t,x)\in \supp(\th_\ep^\NN)\times B$
such that the $F$-orbit
of $(\un t,x)$ has the same $(\un t,x)$ as an accumulation
point. We see that there must exist some $p\in\NN$ such
that $f^p_{\un t}(x)\in B\subset C_i$.
In view of the independence of the permutation on the
choice of $\un t$, we conclude that $C_i$ is sent inside
itself by $f_{\un t}^p$ for all $\un
t\in\supp(\th_\ep^\NN)$. \cqd

It is clear that the cycles obtained above are invariant
sets.
We are now ready to decompose $\mu$ into some simpler
measures. For that we need the following result.

\cle \label{l.restrict} The normalized restriction of a
stationary  measure to an invariant set is a
stationary measure.
 \fle
 \dem See \cite[Lemma~8.2]{Ar1}. \cqd

We define an {\em invariant domain} in $M$ as a finite
collection $(U_0,\dots,U_{p-1})$ of pairwise separated open
sets, that is, $\ov{U}_i\cap\ov{U}_j = \emptyset$ if $i\neq
j$, such that
 $f^k_{\un t}(U_i)\subset U_{(k+i)\bmod p}$
 for all $k\ge1$,
$i=0,\dots,p-1$ and $\un t\in\supp(\th_\ep^\NN)$.

In order to get the  separation  of the connected
components in a cycle, we may unite those components $C_i$
and $C_j$ such that $\ov{C}_i\cap\ov{C}_j\neq\emptyset$ and
observe that the permutation now induced in the new sets by
$f_t$ also does not depend on the choice of
$t\in\supp(\th_\ep)$. In this manner we construct invariant
domains inside the support of any stationary probability
measure.

The next step is to look for \emph{minimal invariant
domains} with respect to the natural order relation of
inclusion of sets.
Let $D=(U_0,\dots,U_{p-1})$ and $D'=(W_0,\dots,W_{q-1})$ be
invariant domains.
On the one hand, $D=D'$ if there are $i,j\in\NN$ such that
$U_{(i+k) \bmod p} = W_{(j+k) \bmod q}$ for all $k\ge1$,
which implies $p=q$ because the open sets that form each
invariant domain are pairwise disjoint. On the other hand,
we say $D\prec D'$ if there are $i,j\in\NN$ such that
$U_{i\bmod p} \subsetneq W_{j\bmod q}$
and $U_{(i+k) \bmod p} \subset W_{(j+k)
\bmod q}$ for all $k\ge1$.

\cle \label{l.minimal}
 In the partially ordered family of
all invariant domains in $M$, with respect to the relation
$\prec$, the number of $\prec$-minimal domains is finite.
Moreover, every invariant domain contains at least one
minimal domain.
 \fle

\dem The proof relies in showing that Zorn's Lemma can be
applied to this partially ordered set and that minimal
domains are pairwise separated. See~\cite[Section~3]{Ar1}.
\cqd

Let us now  fix $x\in M$ and consider
 \begin{equation}
 \label{pushforward1}
 \mu_n(x)=\frac1n\sum_{j=0}^{n-1} (f^j_x)_*\th_\ep^\NN.
 \end{equation}
Since this is a sequence of probability measures on the
compact manifold $M$, then it has weak$^*$ accumulation
points.

\cle \label{p.muinvariant}
 Every weak$^*$ accumulation point
of $\big(\mu_n(x)\big)_n$ is stationary and absolutely
continuous with respect to the Lebesgue measure.
 \fle
\dem Let $\mu$ be a weak$^*$ accumulation point of
$\big(\mu_n(x)\big)_n$. We may write
 $$
 \int \int \vfi \big( f_t(x) \big) \,
d\mu(x) \, d\th_\ep(t) = \int \lim_{k\to+\infty}
\frac1{n_k} \sum_{j=0}^{n_k-1} \int \vfi \left(
f_t\big(f_{\un t}^j(x)\big) \right) \, d\th_\ep^\NN(\un
t)\,d\th_\ep(t)
 $$
 for each continuous $\vfi:M\to\RR$.
Moreover dominated convergence ensures that we may exchange
the limit and the outer integral sign and, by definition of
$f^j_{\un t}(x)$, we get
 $$
 \lim_{k\to\infty} \frac1{n_k}
\sum_{j=0}^{n_k-1} \int \vfi \big( f_{\un t}^{j+1}(x)\big)
 \, d\th_\ep^\NN(\un t) = \int \vfi \, d\mu,
 $$
according to the definition of $\mu$. Thus
(\ref{eq.stationary}) must hold and $\mu$ is stationary.

Noting that $C^0(M,\RR)$ is dense in $L^1(M,\mu)$ with the
$L^1$ norm, we see that (\ref{eq.stationary}) holds for all
$\mu$-integrable functions $\vfi:M\to\RR$. In particular,
if $E\subset M$ is such that $m(E)=0$, then
 \begin{eqnarray*}
 \int 1_E \,d\mu
 &=&\int\int 1_E\big(f_{t}(x)\big) \,d\mu(x)\,d\th_\ep(t) \\
 &=&\int\int 1_E\big(f_{t}(x)\big) \,d\th_\ep(t)\,d\mu(x) \\
 &=&\int\int\int1_E\big(f_{t}(f_{s}(x))\big)\,
       d\th_\ep(t)\,d\mu(x)\,d\th_\ep(s)\\
 &=&\int\int 1_E\big(f^2_{\un t}(x)\big)\,d\th_\ep^\NN(\un t)\,d\mu(x)\\
 &=&\int (f^2_x)_*\th_\ep^\NN(E) \,d\mu(x).
 \end{eqnarray*}
 This process may be iterated to yield
 $$
 \mu(E)=\int (f^{n_0}_x)_*\th_\ep(E)\,d\mu(x)
 $$
and, since $(f^{n_0}_x)_*\th_\ep\ll m$ by nondegeneracy
condition 2, we must have $\mu(E)=0$.
 \cqd

Clearly if $x\in M$ belongs to some set of an
invariant domain $(U_0,\dots,U_{p-1})$,
then $\mu_n(x)$ have supports contained in
$\ov{U}_0\cup\dots\cup\ov{U}_{p-1}$ for all $n\ge1$
and any weak$^*$ accumlation point $\mu$ of
$(\mu_n(x))_n$ is a stationary measure with
$\supp(\mu)\subset\ov{U}_0\cup\dots\cup\ov{U}_{p-1}$.
We will now see these measures are physical.

 \cle
 \label{l.minimal=physical}
If $(U_0,\dots,U_{p-1})$ is a minimal invariant
domain, then there is a unique absolutely
continuous stationary measure $\nu$ such that
$\supp(\nu)\subset\ov{U}_0\cup\dots\cup\ov{U}_{p-1}$.
Moreover, this  $\nu$ is a physical measure
 and $\supp(\nu)=\ov{U}_0\cup\dots\cup\ov{U}_{p-1}$.
 \fle

\dem Let us assume $n_0=1$ for simplicity
(see~\cite[Section~7]{Ar1} for the general case) and let us
consider a stationary absolutely continuous probability
measure $\nu$ with
$\supp(\nu)\subset\ov{U}_0\cup\dots\cup\ov{U}_{p-1}$. We
first show  the ergodicity of $\nu$, in the sense that
$\th_\ep^\NN\times\nu$ is $F$-ergodic. It turns out that to
be $F$-ergodic it suffices that either $\nu(G)=0$  or
$\nu(G)=1$ for every Borel set $G\subset M$
 satisfying
 \begin{equation}
 \label{eq.ergodic}
 1_G(x)=\int 1_G \left(f_t(x)\right)\,d\th_\ep(t)
\end{equation} for $\nu$ almost every $x$ (cf.~\cite{Ar1}
and~\cite{V2}). So let us take $G$ such that $\nu(G)>0$ and
$G$ satisfies the left hand side of~(\ref{eq.ergodic}).
Then it must be $m(G)>0$ because $\nu\ll m$ and there is a
closed set $J\subset G$ such that $m(G\setminus J)=0$ and
also $\nu(G\setminus J)=0$. Hence $J$ also satisfies the
left hand side of~(\ref{eq.ergodic}) because of
nondegeneracy condition 2 (with $n_0=1$), since
 $$
 \int 1_E(f_t(x))\,d\th_\ep(t)=(f_x)_*\th_\ep^\NN(E).
 $$
This means that when $x\in J$ we have $f_t(x)\in J$ for
$\th_\ep$ almost all $t\in\supp(\th_\ep)$. Since a set of
$\th_\ep$ measure 1 is dense in $\supp(\th_\ep)$ (we are
supposing $\th_\ep$ to be positive on open sets) and
$f_t(x)$ varies continuously with $t$,
we see that $f_t(x)\in J$ for all $t\in\supp(\th_\ep)$
because $J$ is closed.
We then  have that the interior of $J$ is nonempty by
condition 1 on random perturbations and we may apply the
methods of decomposition into connected components as
before (Lemma~\ref{pr.cycles}). In this manner we
construct an invariant domain inside $J$ which, in turn, is
inside a minimal invariant domain. This contradicts
minimality and so we conclude that $J$ must contain
$\ov{U}_0\cup\dots\cup\ov{U}_{p-1}$. Thus we have
$\nu(G)=\nu(J)=1$ proving $\th_\ep^\NN\times\nu$ to be
$F$-ergodic.

Now, given $\vfi:M\to\RR$ continuous we consider the map
$\psi=\vfi\circ\pi$ from $T^\NN\times M$ to $\RR$,
where $\pi:T^\NN\times M\to M$ is
the natural projection. The Ergodic Theorem then ensures
 $$
 \lim_{n\to+\infty}\frac1n\sum_{j=0}^{n-1}
 \psi(F^j(\un t,x))=\int \psi \, d(\th_\ep^\NN\times\nu)
 $$
for $\theta_\epsilon^\NN\times\nu$ almost all $ (\un t,x)$,
which is just the same as
 \begin{equation}
 \label{eq.ergodic2}
 \lim_{n\to+\infty}\frac1n\sum_{j=0}^{n-1}
 \vfi(f^j_{\un t}(x))=\int\vfi\,d\nu
 \end{equation}
 for $\theta_\epsilon^\NN\times\nu$ almost all $ (\un
 t,x)$.
 Finally considering the ergodic basin
 $B(\nu)$, defined as the set of points $ x\in M$ for which
 $$
 \lim_{n\to+\infty}
 \frac1n\sum_{j=0}^{n-1}\vfi(f^j_{\un t}(x))=
 \int\vfi\,d\nu
 $$
 for all $\vfi\in C^0(M,\RR)$ and
 $\th_\ep^\NN$ almost every $\un t\in T^\NN$,
it is easy to
see that $B(\nu)$ satisfies 
(\ref{eq.ergodic}) in the place of $G$ and we must have
as before $B(\nu)\supset \ov{U}_0\cup\dots\cup\ov{U}_{p-1}$.

This shows that if another stationary absolutely continuous
probability measure $\tilde \nu$ is such that $\supp(\tilde
\nu)\subset \ov{U}_0\cup\dots\cup\ov{U}_{p-1}$, then the
basins of $\nu$ and $\tilde \nu$ must have nonempty intersection.
Thus these
measures must be equal. Moreover $\nu\big(B(\nu)\big)=1$ and so, by
absolute continuity, $m\big(B(\nu)\big)>0$ and thus $\nu$ is a physical
probability. \cqd


\section{The number of physical measures}
\label{s.physical}

In this section we will prove that the number $l$ of physical
measures is bounded by the number $p$ of SRB measures. Moreover we will present examples of dynamical systems for which $l=p$ and $l<p$.

Let $\mu_1,\ldots,\mu_l$ be the physical measures supported
on the minimal invariant domains  in $M$, which exist by
Lemmas~\ref{pr.cycles} and \ref{l.minimal}
through~\ref{l.minimal=physical}. If $\mu$ is an absolutely
continuous stationary measure, its restrictions to the
minimal invariant domains of $M$, normalized when not equal
to the constant zero measure, are absolutely continuous
stationary measures by Lemma~\ref{l.restrict}. After
Lemma~\ref{l.minimal=physical} these restrictions must be
the physical measures $\mu_1,\ldots,\mu_l$ of the minimal
domains. Hence $\mu$ must decompose into a linear
combination of physical measures.
 Moreover, the union of
$\supp(\mu_1),\dots,\supp(\mu_l)$ must contain
$\supp(\mu)$, except possibly for a $\mu$ null set. In
fact, if the following set function
 $$
\mu-\mu\big(\supp(\mu_1)\big)\mu_1-\cdots
-\mu\big(\supp(\mu_l)\big)\mu_l
 $$
were nonzero, then its normalization $\mu^\prime$ would be
an absolutely continuous stationary measure, and the above
 decomposition could be applied to $\mu^\prime$,
thus giving another minimal domain inside $\supp(\mu)$.
Clearly this cannot happen. We then have a convex linear
decomposition
 \begin{equation}
\label{eq.decomposition}
 \mu=\al_1\mu_1+\dots+\al_l\mu_l
 \end{equation}
 where $\al_i=\mu(\supp(\mu_i))\ge0$ and
$\al_1+\dots+\al_l=1$. We will see that this decomposition
is uniquely defined.

\medskip

We remark that so far we did not use more than the continuity
of the map $f$.
For the next result we assume that $f:M \to M$ is a $C^2$
non-uniformly expanding map whose orbits have slow
approximation to the critical $\cc$ (possibly the emptyset)
with $m(\cc)=0$. This result contains the assertions of the
first two items of Theorem~A (if we think of
$\cc=\emptyset$) and Theorem~C.

\cpr
\label{p.finitemeasures}
 If $\ep>0$ is small enough, then there exist physical measures
 $\mu^\ep_1,\dots,\mu^\ep_l$
(with $l$ not depending on $\ep$) such that
 \begin{enumerate}
 \item for $x\in M$ there is a $\th_\ep^\NN$ mod $0$
 partition $T_1(x),\dots,T_l(x)$ of $T^\NN$ such that
 $$
 \mu_i^\ep=\wlim \frac1n \sum_{j=1}^{n-1} \de_{f_{\un t}^j(x)}
 \quad \mbox{ if and only if } \quad \un t \in T_i(x);
 $$
 \item for each $i=1,\dots,l$ we have
 $$
 \mu_i^\ep=\wlim\frac1n\sum_{j=0}^{n-1}\int (f^j_{\un t})_*
 \big(m\mid B(\mu_i^{\ep})\big) \, d\th_\ep^\NN(\un t),
 $$
where $m\mid B(\mu_i^{\ep})$ is the normalized
restriction of Lebesgue measure to $B(\mu_i^{\ep})$.
 \end{enumerate}
 \fpr
 \dem
Take $x\in M$ and let $\mu$ be a weak$^*$ accumulation
point of the sequence $(\mu_n(x))_n$ defined
in~(\ref{pushforward1}). We will prove that this is  the
only accumulation point of~(\ref{pushforward1}) by showing
that  the values of the $\al_1,\dots,\al_l$ in
decomposition (\ref{eq.decomposition}) depend only on $x$
and not on the subsequence that converges to $\mu$. The
definition of the average in~(\ref{pushforward1}) implies
that there is a subset
 of parameter vectors  $\un t\in\supp(\th_\ep^\NN)$
with positive $\th_\ep^\NN$ measure
for which  there is $j\ge1$ such that $f^j_{\un
t}(x)\in\supp(\mu_i)$.
We define for $i=1,\dots,l$
 $$
 T_i(x)=
 \left\{ \un t\in\supp(\th_\ep^\NN): f^j_{\un t}(x) \in
 \supp(\mu_i) \quad\mbox{for some}\quad j\ge1 \right\}.
 $$
We clearly have
 $$
 T_i(x)=\mbox{$\bigcup$}_{j\ge1} T_i^j(x) \quad\mbox{where}\quad T_i^j(x)=
 \{ \un t\in\supp(\th_\ep^\NN): f^j_{\un t}(x)\in \supp(\mu_i) \}
 $$
 and $T_i^j(x)\subset T_i^{j+1}(x)$ for all $i,j\ge1$, since the supports of
stationary measures are themselves invariant. In addition,
since $\mu$ is a regular (Borel) probability measure, we
may find for each $\eta>0$ an open set $U$ and a closed set $K$
such that $K\subset\supp(\mu_i)\subset U$ with
$\mu(U\setminus K)<\eta$ and $\mu(\partial U)=\mu(\partial
K)=0$. In fact, there is an at most countable number of
$\delta$-neighborhoods of $\supp(\mu_i)$ whose boundaries have
positive $\mu$ measure, and likewise for the compacts coinciding with the complement of the $\delta$-neighborhood of $M\setminus\supp(\mu_i)$. Then,
taking $\alpha_i=\mu(\supp(\mu_i))$ we have
 \begin{eqnarray*}
\al_i+\eta\ge \mu(U) &=&
\lim_{k\to+\infty}\frac1{n_k}\sum_{j=0}^{n_k-1}
\th_\ep^\NN\{ \un t\in T^\NN: f^j_{\un t}(x)\in U \} \\ &
\ge & \limsup_{k\to+\infty} \frac1{n_k}\sum_{j=0}^{n_k-1}
\th_\ep^\NN\big(T_i^j(x)\big) \end{eqnarray*} for some sequence of
integers $n_1<n_2<n_3<\cdots$, and likewise for
\begin{eqnarray*} \al_i-\eta\le \mu(K)
 &=&
 \lim_{k\to+\infty}\frac1{n_k}\sum_{j=0}^{n_k-1}
 \th_\ep^\NN\{ \un t\in T^\NN:
 f^j_{\un t}(x)\in K \}
\\
 & \le &
 \liminf_{k\to+\infty}
 \frac1{n_k}\sum_{j=0}^{n_k-1}
 \th_\ep^\NN\big(T_i^j(x)\big),
\end{eqnarray*}
where $\eta>0$ is arbitrary.
This shows
 $$
 \al_i=\mu(\supp(\mu_i))= \lim_{k\to\infty}
 \frac1{n_k}\sum_{j=0}^{n_k-1} \th_\ep^\NN\big(T_i^j(x)\big).
 $$
We also have
 $$
\th_\ep^\NN\big(T_i(x)\big)=\lim_{j\to\infty}\th_\ep^\NN\big(T_i^j(x)\big)
=\lim_{n\to\infty}\frac1n\sum_{j=0}^{n-1}
\th_\ep^\NN\big(T_i^j(x)\big) =\al_i
 $$
which shows that the $\al_i$ depend only on the random
orbits of $x$ and not on the particular sequence $(n_k)_k$.
Thus we see that the sequence of measures in
(\ref{pushforward1}) converges in the weak$^*$ topology.
Moreover the sets $T_1(x),\dots,T_l(x)$ are pairwise disjoint by
definition and their total $\th_\ep^\NN$ measure equals
$\al_1+\cdots+\al_l=1$, thus forming a $\th_\ep^\NN$ modulo
zero partition of $T^\NN$. We observe that if $\un t\in
T_i(x)$, then $f^n_{\un t}(x)\in \supp(\mu_i)\subset
B(\mu_i)$ for some $n\ge1$ and $i=1,\dots,l$. This means
this $\th_\ep^\NN$ modulo zero partition of $T^\NN$
satisfies the first item  of the proposition.

Now fixing $i=1,\dots,l$, for all $x\in B(\mu_i)$ (the
ergodic basin of $\mu_i$) it holds that $$
\lim_{n\to+\infty}\frac1n\sum_{i=0}^{n-1} \vfi(f_{\un
t}^j(x)) = \int \vfi\, d\mu_i $$ for $\th_{\ep}^\NN$ almost
every $\un t\in T^\NN$. Recall that $m(B(\mu_i))>0$ by the
definition of physical measure. Using dominated convergence
and integrating both sides of the above equality twice,
first with respect to the Lebesgue measure $m$, and then
with respect to $\th_{\ep}^\NN$, we arrive at the statement
of item 2.

Recall that up until now the noise level $\ep>0$ was kept
fixed. For small enough $\ep>0$ the measures
$\mu_i=\mu_i^\ep$ depend on the noise level, but we will
see that the number of physical measures is constant.

Fixing $i\in\{1,\dots,l\}$ we let $x$ in the interior of
$\supp(\mu_i^\ep)$ be such that the orbit $(f^j(x))_j$ has
infinitely many hyperbolic times. Recall that $f\equiv
f_{t^*}$ is non-uniformly expanding (possibly with
criticalities). Then there is a big enough hyperbolic time
$n$ so that $V_n(\un t^*, x)\subset\supp(\mu_i^\ep)$, by
Proposition~\ref{p.contr}, where we take $\un
t^*=(t^*,t^*,t^*,\dots)$. Since $t^*\in\supp(\th_\ep)$ and
$\supp(\mu_i^\ep)$ is invariant under $f_t$ for all
$t\in\supp(\th_\ep)$, we must have
 $$
 f^n_{\un t^*}\big(V_n(\un
t^*,x)\big)= B\big(f^n_{t^*}(x),\de_1\big)
\subset\supp(\mu_i^\ep),
 $$
where $\de_1>0$ is the constant given by
Proposition~\ref{p.contr} and
$B\big(f^n_{t^*}(x),\de_1\big)$ is the ball of radius
$\de_1$ around $f^n_{t^*}(x)$.

 On the
one hand, we deduce that the number $l=l(\ep)$ is bounded
from above by some uniform constant $N$ since $M$ is
compact. On the other hand, since each  invariant set  must
contain some physical measure (by Lemma~\ref{l.minimal}),
we see that for $0<\ep^\prime<\ep$ there must be some
physical measure $\mu^{\ep^\prime}$ with
$\supp(\mu^{\ep^\prime})\subset\supp(\mu^\ep)$. In fact
$\supp(\mu^\ep)$ is invariant under $f_t$ for every
$t\in\supp(\th_{\ep^\prime})\subset\supp(\th_\ep)$. This
means the number $l(\ep)$ of physical measures is a
nonincreasing function of $\ep>0$. Thus we conclude that
there must be $\ep_0>0$ such that $l=l(\ep)$ is constant
for $0<\ep<\ep_0$, ending the proof of the proposition.
\cqd

\cre \label{r.transitivity} Observe that if the map $f:M\to
M$ is transitive, then every stationary measure must be
supported on the whole of $M$, since the support is
invariant and has nonempty interior. According to the
discussion above, there must be only one such stationary
measure, which must be physical. \fre

We note that the number $l$ of physical measures for small
$\ep>0$ and the number $p$ of SRB measures for $f$ are
obtained by different existential arguments. It is natural
to ask if there is any relation between $l$ and $p$.

\cpr  If $p\ge1$ is the number of SRB measures of $f$ and
$l\ge1$ is the number of physical measures of the random
perturbation of $f$, then for $\ep>0$ small enough we have $l\le
p$.
 \fpr
 \dem
We start by observing that if $p=1$ then every weak$^*$
accumulation point of a family $(\mu_i^\ep)_{\ep>0}$ of
physical measures when $\ep\to0$ must equal the unique SRB
measure $\mu_1$ for $f$. Hence the weak$^*$ limit of
$(\mu_i^\ep)_{\ep>0}$ when $\ep\to0$ exists and equals
$\mu_1$ for all $i=1,\ldots,l$. Then there must be a single
physical measure $\mu_1^\ep$ for all small enough $\ep>0$.
In fact, let us assume there are distinct families
$(\mu_1^\ep)_{\ep>0}$ and $(\mu_2^\ep)_{\ep>0}$ of physical
measures as given by Proposition~\ref{p.finitemeasures}.
Take $x\in\supp(\mu_1)$ and a sequence of continuous maps
$\vfi_n:M\to\RR$ such that
 $$ \vfi_n\ge0\quad,\quad
\vfi_n\mid B(x, 1/n)\equiv 1 \qand \vfi_n \mid \big(
M\setminus B(x,2/n)\big) \equiv 0,
 $$ for all big $n\in\NN$, where $B(x,r)$ denotes the ball of radius
 $r$ around $x$ for each $r\ge0$. If we fix $n$, then $\mu_1(\vfi_n)>0$. Thus for all
small enough $\ep>0$ we must have $\mu_1^\ep(\vfi_n)>0$ and
$\mu_2^\ep(\vfi_n)>0$. Hence the distance between
$\supp(\mu_1^\ep)$ and $\supp(\mu_2^\ep)$ is smaller than
$2/n$.\ Since $n$ may be arbitrarily large we see that
$\supp(\mu_1^\ep)$ and $\supp(\mu_2^\ep)$ get arbitrarily
close when $\ep\to0$. Thus $\mu_1^\ep$ and $\mu_2^\ep$ must
coincide for small $\ep>0$ because
Proposition~\ref{p.finitemeasures} and its proof show that
these families, when distinct, are at a uniform distance
apart (since $\{\supp(\mu_1^\ep)\}_{\ep>0}$ and
$\{\supp(\mu_2^\ep)\}_{\ep>0}$ are nested families of
pairwise disjoint compact sets).

In general, if $p>1$ then the weak$^*$ accumulation points
of a family $(\mu_i^\ep)_{\ep>0}$ for $i=1,\ldots, l$
when $\ep\to0$ are convex linear combinations
$\al_1\mu_1+\cdots+\al_p\mu_p$ of the $p$ SRB measures of
$f$. The preceding argument implies that two distinct
families $(\mu_i^\ep)_{\ep>0}$ and $(\mu_j^\ep)_{\ep>0}$ of
physical measures cannot have weak$^*$ accumulation points
expressed as linear convex combinations $$
\al_1\mu_1+\cdots+\al_p\mu_p\qand
\al_1^\prime\mu_1+\cdots+\al_p^\prime\mu_p $$ with both
$\al_k$ and $\al_k^\prime$ nonzero for some $k=1,\ldots,p$
(take $x\in\supp(\mu_k)$ and repeat the arguments in the
above paragraph). Hence $l\le p$. \cqd

The reverse inequality does not hold in general, as the
following examples show: it is possible for two distinct
SRB measures to have intersecting supports and, in this
circumstance, the random perturbations will mix their
basins and there will be some physical measure whose
support overlaps the supports of both SRB measures.

\begin{figure}[htb]
\begin{center}
\includegraphics[width=8cm]{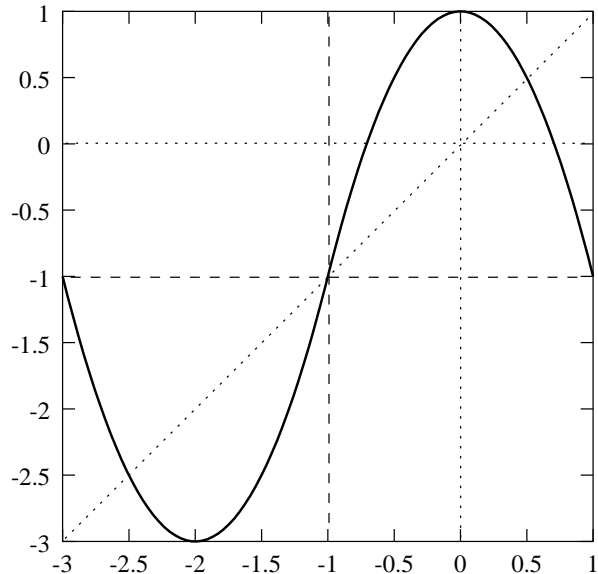}
\caption{\label{fig.1}
            map for which $1=l<p=2$}
\end{center}
\end{figure}

\medskip

The first example is the map $f:[-3,1]\to[-3,1]$ whose
graph is figure~\ref{fig.1}: $$ f(x)=  \left\{
\begin{array}{lcl} 1-2x^2 & \mbox{if} & -1\le x \le 1 \\
2(x+2)^2 -3 & \mbox{if} & -3 \le x \le -1 \end{array}
\right. . $$ The dynamics of $f$  on $[-1,1]$ and $[-3,-1]$
 is conjugated to the tent map $T(x)=
1-2|x|$ on $[-1,1]$. Thus understanding $f$ as a circle map
through the identification $S^1=[-3,1]/\{-3,1\}$, this is a
non-uniformly expanding map with a critical set satisfying
conditions (S1)-(S3) and
 there are two ergodic absolutely continuous (thus SRB)
invariant measures $\mu_1, \mu_2$ whose supports are
$[-3,-1]$ and $[-1,1]$ respectively. Moreover defining
 $\Phi(t)= R_t\circ f,
 $
 where $R_t:S^1\to S^1$ is the
rotation of angle $t$ and $\th_\ep=(2\ep)^{-1}(m\mid
[-\ep,\ep])$ for small $\ep>0$, we have that $\{ \Phi,
(\th_\ep)_{\ep>0} \}$ is a random perturbation satisfying
nondegeneracy conditions 1 and 2. Since
$\supp(\mu_1)\cap\supp(\mu_2)=\{-1\}$ we have that for
$\ep>0$  small enough there must be a single physical
measure $\mu^\ep$. Indeed, by property (P) any weak$^*$
accumulation point of a family of physical measures must
have $-1$ in its support.
 \begin{figure}[htb]
\begin{center} \includegraphics[width=8cm]{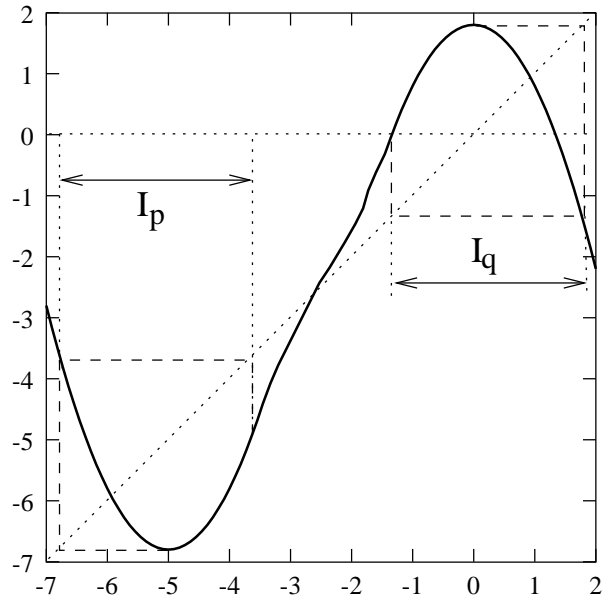}
\caption{\label{fig.2} map for which $l=p=2$} \end{center}
\end{figure}

\medskip

The second example is defined on the interval $I=[-7,2]$.
We take the map $q_a(x)=a-x^2$ on $[-2,2]$ for some
parameter $a\in(1,2)$  satisfying Benedicks-Carleson
conditions (see~\cite{BC1} and~\cite{BC2}), and  the
``same" map on $[-7,-3]$ conveniently conjugated:
$p_a(x)=(x+5)^2-5-a$. Then the two pieces of graph are
glued together in  such a way that we obtain a smooth map
$f:I\to I$ sending $I$ into its interior, as
figure~\ref{fig.2} shows. The intervals
$I_q=[q_a^2(0),q_a(0)]$ and $I_p=[p_a(-5),p_a^2(-5)]$ are
forward invariant for $f$, and then we can find slightly
larger intervals $I_1\supset I_p$ and $I_2\supset I_q$
that become trapping regions for $f$. So, taking
 $ \Phi(t)= f+t, $
 and $\theta_\ep$ as in the previous
example with $0<\ep<\ep_0$ for some $\ep_0>0$ small enough,
then $\{ \Phi, (\th_\ep)_{\ep}\}$ is a random perturbation
of $f$ leaving the intervals  $I_1$ and $I_2$ invariant by
each $\Phi(t)$. Moreover, Lebesgue almost every $x\in I$
eventually arrives at one of these intervals. Then
by~\cite{BC1} and \cite{BY} the map $f$ is non-uniformly
expanding and has two SRB measures with supports contained
in each trapping region. Finally $f$ admits two distinct
physical measures whose supports are contained  in $I_1$
and $I_2$ respectively, for $\ep_0>0$ small enough, see
\cite{BaV}.



\section{Stochastic stability} \label{s.stochastic}

In this section we will prove the first item of
Theorem~\ref{t.stc} and Theorem~\ref{t.sts}. The second
item of Theorem~\ref{t.stc} may be obtained in the same way
as Theorem \ref{t.sts}, if we think of $\cc$ as being equal
to the empty set and take into account Remark~\ref{r.hyp}.

We start by proving the first item of Theorem \ref{t.stc}.
Assume that $f$ is a stochastically stable non-uniformly
expanding  local diffeomorphism. We know from
Proposition~\ref{p.finitemeasures} that there is a finite
number of physical measures $\mu_1^\ep,\dots \mu_l^\ep$ and
for each $x\in M$ there is a $\th_\ep^\NN$ mod $0$
partition $T_1(x),\dots,T_l(x)$ of $T^\NN$ for which
 $$
 \mu_i^\ep=\wlim \frac1n \sum_{j=1}^{n-1}\de_{f_{\un t}^j(x)}
 \quad \mbox{for each}\quad \un t \in T_i(x).
 $$
Furthermore, since we are taking $f$ a local
diffeomorphism, then $\log\|(Df)^{-1}\|$ is a continuous
map. Thus, we have for each $x\in M$ and $\th^\NN_\ep$
almost every $\un t\in T^\NN$
 $$
 \lim_{n\ra \infty}\frac{1}{n}\sum_{j=0}^{n-1}
 \log\|Df\big(f_{\un t}^j(x)\big)^{-1}\|=
 \int\log\|(Df)^{-1}\|d\mu_i^\ep
 $$
for some physical measure $\mu_i^\ep$ with $1\le i\le l$.
Hence, for proving the non-uniform expansion of $f$ on
random orbits it suffices to show that there is $c_0>0$
such that if $\mu^\ep=\mu_i^\ep$ for some $1\le i\le l$
then
 $$
 \int\log\|(Df)^{-1}\|d\mu^\ep<c_0\quad\mbox{for small $\ep>0$}.
 $$

 \cle\label{l.aprox}
Let $\vfi\colon M\ra \RR$ be a continuous map.
Given $\de>0$ there is $\ep_0>0$ such that if $\ep\leq\ep_0$, then
 $$
 \left|\int \vfi d\mu^\ep -
 \int \vfi d\mu_\ep\right|<\de,
 $$
for some absolutely continuous $f$-invariant probability
measure $\mu_\ep$.
 \fle
 \dem
We will use the following auxiliary result: \emph{Let $X$
be a compact metric space, $K\subset X$ a closed (compact)
subset and $(x_t)_{t>0}$ a curve in $X$ (not necessarily
continuous) such that all its accumulation points (as $t\to
0^{+}$) lie in $K$. Then for every open neighborhood $U$ of
$K$ there is $t_0>0$ such that $x_t\in U$ for every
$0<t<t_0$. } Indeed, supposing not, there is a sequence
$(t_n)_{n}$ with $t_n\to 0^{+}$ when $n\to\infty$ such that
$x_{t_n}\not\in U$. Since $X$ is compact this means that
$(x_t)_{t>0}$ has some accumulation point in $X\setminus
U$, thus outside $K$, contrary to the assumption.

Now, the space $X=\PP(M)$ of all probability measures in
$M$ is a compact metric space with the weak$^{*}$ topology,
and the convex hull $K$ of the (finitely many) SRB measures
of $f$ is closed. Hence, considering the curve
$(\mu_\ep)_\ep$ in $\PP(M)$, we are in the context of the
above result, since we are supposing $f$ to be
stochastically stable. A metric on $X$ topologically
equivalent to the weak$^{*}$ topology may be given by
 $$
\distp(\mu,\nu)=\sum_{k=1}^\infty \frac1{2^n} \left| \int
\varphi_n \, d\mu - \int \varphi_n \, d\nu \right|
 $$
 where
$\mu,\nu\in\PP(M)$ and $(\varphi_n)_{n\ge 1}$ is a
dense sequence of functions in $C^0(M,\RR)$,
see~\cite{Man}.

Let $\varphi:M\to\RR$ continuous be given and let us fix
some $\de>0$. There must be $n\in\NN$ such that $\| \varphi
- \varphi_n \|_0 < \de/3$ and, by the auxiliary result in
the beginning of the proof,  there exists, for some $\ep_0>0$ and
every $0<\ep<\ep_0$, a probability measure
$\mu_\ep\in\PP(M)$ for which $\distp(\mu^\ep,\mu_\ep)<\de
(3\cdot 2^n)^{-1}$. This in particular means that
 $$
\frac1{2^n} \left| \int \varphi_n \, d\mu^\ep - \int
\varphi_n \, d\mu_\ep \right| < \frac{\de}{3\cdot 2^n},
 $$
by the
definition of the distance $\distp$, which implies
 $$
 \left| \int \varphi_n \, d\mu^\ep - \int
\varphi_n \, d\mu_\ep \right| < \frac{\de}3.
 $$
Hence we get
 \begin{eqnarray*}
\lefteqn{\left| \int \varphi \, d\mu^\ep - \int \varphi \,
d\mu_\ep \right|\leq}\hspace{.5cm}\\
 & \le &
 \left| \int \varphi \, d\mu^\ep - \int
\varphi_n \, d\mu^\ep \right| + \left| \int \varphi_n \,
d\mu^\ep - \int \varphi_n \, d\mu_\ep \right| + \left| \int
\varphi_n \, d\mu_\ep - \int \varphi \, d\mu_\ep \right| \\
 & < & \frac{\de}3 +\frac{\de}3 +\frac{\de}3 =
\de, \end{eqnarray*}
which completes the proof of the lemma.\cqd

Now we take $\vfi=\log\|(Df)^{-1}\|$ and $\de=c/2$ in the
previous lemma, where $c>0$ is the constant given by the
non-uniform expansion  of $f$ (recall (\ref{liminf1})). For
each $\ep\leq\ep_0$ let $\mu_\ep$ be the measure given by
Lemma~\ref{l.aprox}. Since property (P) holds, there are
real numbers $w_1(\ep),\dots ,w_p(\ep)\geq 0$ with
$w_1(\ep)+\cdots+w_p(\ep)=1$ for which
$\mu_\ep=w_1(\ep)\mu_1+\cdots+w_p(\ep)\mu_p$. Since each
$\mu_i$ is an SRB measure for $1\leq i\leq p$, we have for
Lebesgue almost every $x\in B(\mu_i)$
 $$
\int\log\|(Df)^{-1}\|d\mu_i=\lim_{n\ra +\infty}\frac{1}{n}
\sum_{j=0}^{n-1}\log\|Df(f^j(x))^{-1}\|\leq -c<0.
 $$
 This implies
 $$
 \int\log\|(Df)^{-1}\|d\mu_\ep\leq -c,
 $$
and so, by Lemma \ref{l.aprox} and the choice of $\delta$,
 $$
 \int\log\|(Df)^{-1}\|d\mu^\ep\leq -c/2.
 $$
This completes the proof of the first item of Theorem B.

\medskip

Now we go into the proof of Theorem \ref{t.sts}.
In order to prove that $f$ is stochastically stable, and
taking into account property (P), it suffices to prove
that the weak$^*$ accumulation points of any family
$(\mu^\ep)_{\ep>0}$, where each $\mu^\ep$ is a physical
measure of level $\ep$, are absolutely continuous with respect
to the Lebesgue measure. Let $\mu^\ep$ be a physical measure
of level $\ep$ for some small $\ep>0$
 and define
for each $n\geq 1$
 $$ \mu_n^\ep = \frac1n\sum_{j=0}^{n-1}\int (f^j_{\un t})_{*}
 \big(m\mid{B(\mu^\ep)\big)}\, d\th_\ep^\NN(\un t).
 $$
We know from Proposition~\ref{p.finitemeasures} that each
$\mu^\ep$ is the weak$^*$ limit of the sequence
$(\mu_n^\ep)_n$. We will prove Theorem~\ref{t.sts} by
providing some useful estimates on the densities of the
measures $\mu_n^\ep$. Define for each $\un t\in T^\NN$ and
$n\geq 1$
 $$
 H_n(\un t)=\{ x\in  B(\mu^\ep)\colon \mbox{ $n$ is a
 $(\al,\delta)$-hyperbolic time for $(\un t,x)$ }\},
 $$
 and
 $$
  H^*_n(\un t) =\{ x\in B(\mu^\ep)\colon \mbox{ $n$ is the first
 $(\al,\delta)$-hyperbolic time for $(\un t,x)$ }\}.
 $$
 $ H^*_n(\un t)$ is precisely the set of those points
 $x\in B(\mu^\ep)$ for which $h_\ep(\un t,x)=n$
 (recall the definition of
the map $h_\ep$).
 For
$n,k\geq 1$ we also define $R_{n,k}(\un t)$ as the set of
those points $x\in M$ for which $n$ is a
$(\al,\delta)$-hyperbolic time and $n+k$ is the first
$(\al,\delta)$-hyperbolic time  after $n$, i.e. $$
R_{n,k}(\un t)= \left\{x\in H_n(\un t)\colon \:f^n_{\un
t}(x)\in  H^*_k(\si^n\un t)\: \right\}, $$ where $\si\colon T^\NN \rightarrow  T^\NN$ is the shift map
$\si(t_1,t_2,\dots)=(t_2,t_3,\dots)$.
  Considering the
measures
 $$
 \nu^\ep_n=\int ( f_{\un
t}^n)_*\big(m\mid H_n(\un t)\big)d\th_\ep^\NN(\un t)
 $$
 and
 $$
\eta_n^\ep=\sum_{k=2}^\infty\sum_{j=1}^{k-1}\int (f_{\un
t}^{n+j})_*\big(m\mid R_{n,k}(\un t)\big) d\th_\ep^\NN(\un
t),
 $$
we may write
 $$
 \mu_n^\ep\leq
\frac{1}{n}\sum_{j=0}^{n-1}(\nu_j^\ep+\eta_j^\ep).
 $$

\cpr\label{p.dens}
 There is a constant $C_2>0$ such that for every $n\geq 0$
 and $\un t\in T^\NN$
 $$
 \frac{d}{dm}(f_{\un t}^n)_*\big(m\mid H_n(\un t)\big)\leq C_2.
 $$
 \fpr
 \dem Take $\de_1>0$ given by Proposition \ref{p.contr}. It is
sufficient to prove that there is some uniform constant
$C>0$ such that if $A$ is a Borel set in $M$ with diameter
smaller than $\delta_1/2$ then
 $$
 m\big(f_{\un t}^{-n}(A)\cap  H_n(\un t)\big)\leq C m(A).
 $$
 Let $A$ be a Borel set in $M$ with diameter smaller than
$\delta_1/2$ and $B$ an open  ball of radius $\delta_1/2$
containing $A$. We may write
 $$
 f_{\un t}^{-n}(B)=\bigcup_{k\geq 1}B_k,
 $$
where $(B_k)_{k\geq 1}$ is a (possibly finite) family of
two-by-two disjoint open sets in $M$. Discarding those
$B_k$ that do not intersect $H_n(\un t)$, we choose for
each $k\geq 1$ a point $x_k\in H_n(\un t)\cap B_k$. For
$k\geq 1$ let $V_n(\un t,x_k)$ be the neighborhood of $x_k$
in $M$ given by Proposition \ref{p.contr}. Since $B$ is
contained in $B\big(f_{\un t}^n(x_k), {\de_1}\big)$, the
ball of radius $\de_1$ around $f_{\un t}^n(x_k)$, and
$f_{\un t}^n$ is a diffeomorphism from $V_n(\un t,x_k)$
onto $B\big(f_{\un t}^n(x_k), {\de_1}\big)$, we must have
$B_k\subset V_n(\un t,x_k)$ (recall that by our choice of
$B_k$ we have $f_{\un t}^n(B_k)\subset B$). As a
consequence of this and Corollary \ref{c.dist}, we have for
every $k$ that the map $f_{\un t}^n\mid B_k\colon
B_k\rightarrow B$ is a diffeomorphism with bounded
distortion:
 $$ \frac{1}{C_1}\leq \frac{|\det Df_{\un
t}^n(y)|}{|\det Df_{\un t}^n(z)|} \leq C_1 $$ for all
$y,z\in B_k$.
 This finally
gives \begin{eqnarray*} m\big(f_{\un t}^{-n}(A)\cap H_n(\un
t)\big) &\leq & \sum_{k}m\big(f_{\un t}^{-n}(A\cap B)\cap
B_k\big)\\ &\leq & \sum_{k}C_1\frac{m(A\cap
B)}{m(B)}m(B_k)\\ &\leq & C_2 m(A), \end{eqnarray*} where
$C_2>0$ is a constant only depending on $C_1$, on the
volume of the ball $B$ of radius $\delta_1/2$, and on the
volume of $M$.\cqd

It follows from Proposition \ref{p.dens} that
 \begin{equation}\label{dens}
\frac{d\nu_n^\ep}{dm}\leq C_2
 \end{equation}
 for every $n\geq 0$ and small $\epsilon>0$. Our goal now is to
 control the density of the measures $\eta_n^\ep$
 in such a way that we may assure the absolute
 continuity of the weak$^*$ accumulation points
 of the measures $\mu^\ep$ when $\ep$ goes to zero.

\cpr\label{p.dens2} Given $\zeta>0$, there is
$C_3(\zeta)>0$ such that for  every $n\geq 0$ and $\ep>0$
we may bound $\eta_n^{\ep}$ by the sum of two non-negative
measures, $\eta_n^{\ep} \leq \omega^{\:\ep}+\rho^{\:\ep}$,
with
 $$
 \frac{d\omega^{\:\ep}}{dm}\leq C_3(\zeta)\qand
 \rho^{\:\ep}(M)<\zeta.
 $$
 \fpr
 \dem
 Let $A$ be some Borel set in $M$.
We have for each $n\geq 0$
 \begin{eqnarray*}
 \eta_n^\ep(A)
 &=&
 \sum_{k=2}^\infty\sum_{j=1}^{k-1} \int
 m\big( f_{\un t}^{-n-j}(A)\cap R_{n,k}(\un t)\big)
 d\th_\ep^\NN(\un t)\\
 &\leq &
 \sum_{k=2}^\infty\sum_{j=1}^{k-1} \int m\big( f_{\un
 t}^{-n} \big(f_{\si^n\un t}^{-j}(A)\cap H^*_k(\si^n\un t)\big)\cap
 H_n(\un t)\big)d\th_\ep^\NN(\un t)\\
 &=&
 \sum_{k=2}^\infty\sum_{j=1}^{k-1} \int\nu_n^\ep\big(
 f_{\si^n\un t}^{-j}(A)\cap H^*_k(\si^n\un t)\big) d\th_\ep^\NN(\un
 t) \\
 &\leq &
 \sum_{k=2}^\infty\sum_{j=1}^{k-1} C_2\int
 m\big( f_{\un t}^{-j}(A)\cap
 H^*_k(\un t)\big)d\th_\ep^\NN(\un t).
 \end{eqnarray*}
 (in this last inequality we used that $\th^\NN_\ep$
is $\si$-invariant and estimate (\ref{dens}) above). Let
now $\zeta>0$ be some fixed small number. Since we are
assuming $(h_\ep)_\ep$ with uniform  $L^1$-tail this means
that there is some integer $N=N(\zeta)$ for which
 $$
 \sum_{j=N}^{\infty}k\int
 m\big(H^*_k(\un t)\big)d\th^\NN_\ep(\un t)<\frac{\zeta}{C_2}.
 $$
We take
 $$
 \omega^{\:\ep}=C_2\sum_{k=2}^{N-1}\sum_{j=1}^{k-1}
 \int(f_{\un t}^j)_*\big(m\mid H^*_k(\un t)\big)d\th_\ep^\NN(\un t)
 $$
 and
 $$
 \rho^{\:\ep}=C_2\sum_{k=N}^\infty\sum_{j=1}^{k-1}
 \int(f_{\un t}^j)_*\big(m\mid H^*_k(\un t)\big)d\th_\ep^\NN(\un t).
 $$
For this last measure we have
 $$
 \rho^{\:\ep}(M)=C_2\sum_{k=N}^\infty\sum_{j=1}^{k-1} \int
 m\big( H^*_k(\un t)\big)d\th_\ep^\NN(\un t) \leq C_2
 \sum_{k=N}^\infty k\int m\big( H^*_k(\un t)\big)d\th_\ep^\NN(\un t)
 <\zeta.
 $$
On the other hand, it follows from the definition of
$(\al,\de)$-hyperbolic times that there is some constant
$a=a(N)>0$ such that
 $\dist\big(H_k(\un t), \cc\big)\geq a
$ for $1\leq k\leq N$. Defining $\Delta\subset M$ as the
set of those points in $M$ whose distance to $\cc$ is
greater than $a$, we have
 $$
 \omega^{\:\ep}\leq C_2\sum_{k=2}^{N-1}\sum_{j=1}^{k-1}
 \int(f_{\un t}^j)_*(m\mid \Delta)\:d\th_\ep^\NN(\un t),
 $$
and this last measure has density bounded by some uniform
constant, as long as we take the maps $f_t$ in a
sufficiently small neighborhood of $f$ in the $C^1$
topology.
 \cqd

 It follows from this last proposition and (\ref{dens}) that the weak$^*$
 accumulation points of $\mu^\ep$ when $\ep\to0$ cannot
 have singular part, thus being absolutely continuous with
 respect to the Lebesgue measure. Moreover, the
  weak$^*$ accumulation points of a family of
 stationary measures are always $f$-invariant
  measures, cf. Remark~\ref{re.accinvariant}. This together with (P)
 gives the stochastic stability of $f$.


\section{Applications}\label{s.applications}

In this section we will apply  Theorems \ref{t.stc} and
\ref{t.sts} to certain examples of
 non-uniformly expanding maps.
Before we explicit the examples we have in mind let us give
a practical criterion for proving that the family of
hyperbolic time maps $(h_\ep)_\ep$ has uniform $L^1$-tail.

If we look at the proof of Proposition \ref{p.hyp} we
realize that what we did was fixing some positive number
$c_0$ smaller than $c$, and then, for $\th_\ep^\NN\times m$
almost every $(\un t,x)\in T^\NN\times M$, we took a
positive integer $N_\ep=N_\ep(\un t,x)$ for which $$
\sum_{j=0}^{N_\ep-1}\log\|Df(f_{\un t}^j(x))^{-1}\|\leq
-c_0N_\ep \quad \mbox{and} \quad
\sum_{j=0}^{N_\ep-1}-\log\dist_\delta(f_{\un t}^j(x),
\cc)\leq \gamma N_\ep, $$ for suitable choices of
$\delta>0$ and $\gamma>0$. This permits us to introduce a
map
 $$N_\ep\colon T^\NN\times M\ra \ZZ^+$$
whose existence provides a first hyperbolic time map
$$
h_\ep\colon T^\NN\times M\ra \ZZ^+\quad\mbox{with}\quad h_\ep \leq N_\ep
$$
(recall the proof of Proposition \ref{p.hyp}).
Thus, the integrability of the map $h_\ep$ is implied by the
integrability of the map $N_\ep$, which is in practice
easier to handle.

In the examples we are going to study below we will show
that there is a sequence of positive real numbers
$(a_k^\ep)_k$ for which
 $$ (\th_\ep^\NN\times
m)\left(\big\{(\un t,x)\in T^\NN\times M\colon N_\ep(\un
t,x)>k\big\}\right)\leq a_k^\ep \quad \mbox{and} \quad
\sum_{k=1}^{\infty} ka_k^\ep <\infty,
 $$
This gives the integrability of $h_\ep$ with respect to the
measure $\th_\ep^\NN\times m$. The fact the family
$(h_\ep)_\ep$ has uniform  $L^1$-tail can be proved by
showing that the sequence $(a_k^\ep)_k$ may be chosen not
depending on $\ep>0$.

Now we are ready for the applications of Theorems
\ref{t.stc} and \ref{t.sts}. We will describe firstly a
class of local diffeomorphisms introduced in \cite[Appendix
A]{ABV} that satisfies the hypotheses of Theorem
\ref{t.stc}, and then a class of maps (with critical sets)
introduced in \cite{V} satisfying the hypotheses of Theorem
\ref{t.sts}.


\subsection{Local diffeomorphisms} \label{ss.localdiffeos}
Now we follow~\cite[Appendix A]{ABV} and describe robust
classes of maps (open in the $C^2$ topology) that are
non-uniformly expanding local diffeomorphisms and
stochastically stable. Let $M$ be a compact Riemannian
manifold and consider
  \[
  \begin{array}{rccl}
  \Phi:& T &\longrightarrow&  C^2(M,M)\\
  & t &\longmapsto & f_t
  \end{array}
  \]
 a continuous family of $C^2$ maps, where $T$
is a metric space.
 We begin with an essentially
combinatorial lemma.

\cle \label{l.frequency} Let $p,q\ge1$ be integers and
$\si>q$ a real number. Assume $M$ admits a measurable cover
$\{ B_1,\ldots,B_p,B_{p+1},\ldots,B_{p+q}\}$ such that for
all $t\in T$ it holds \begin{enumerate} \item $|\det
Df_t(x)|\ge \si$ for all $x\in B_{p+1}\cup\ldots\cup
B_{p+q}$; \item $\left(f_t\mid B_i\right)$ is injective for
all $i=1,\ldots,p$. \end{enumerate} Then there is $\zeta>0$
such that for every Borel probability $\th$ on $T$ we have
\begin{equation}\label{eq.frequency} \#\{0\le j < n:
f^j_{\un t}(x) \in B_1\cup\ldots\cup B_p \} \ge \zeta  n
\end{equation} for $\th^\NN\times m$ almost all $(\un
t,x)\in T^\NN\times M$ and large enough $n\ge1$. Moreover
the set $I_n$ of points $(\un t,x)\in T^\NN\times M$ whose
orbits do not spend a fraction $\zeta$ of the time in
$B_1\cup\ldots\cup B_p$ up to iterate $n$ is such that
$(\th^\NN\times m)(I_n)\le \tau^n$ for some $0<\tau<1$ and
for large $n\ge1$. \fle

\dem
Let us fix $n\ge1$ and $\un t\in T^\NN$.
For a sequence
$\un i=(i_0,\ldots,i_{n-1})\in\{1,\ldots,p+q\}^n$
we write
$$
[\un i]=
B_{i_0}\cap (f^1_{\un t})^{-1}(B_{i_1})\cap\cdots
\cap (f^{n-1}_{\un t})^{-1}(B_{i_{n-1}})
$$
and define
$g(\un i)=\#\{ 0\le j < n : i_j \le p \}$.

We start by observing that for $\zeta>0$ the number of
sequences $\un i$ such that $g(\un i)<\zeta n$ is bounded
by $$ \sum_{k<\zeta n} {n \choose k} p^k q^{n-k} \le
\sum_{k\le \zeta n} {n \choose k} p^{\zeta n} q^n. $$ Using
Stirling's formula (cf.~\cite[Section 6.3]{BV}) the
expression on the right hand side is bounded by $(e^\gamma
p^\zeta q)^n$, where $\gamma>0$ depends only on $\zeta$ and
$\gamma(\zeta)\to 0$ when $\zeta\to0$.

Assumptions 1 and 2 ensure $m([\un i])\le
\si^{-(1-\zeta)n}$ (recall that $m(M)=1$). Hence the
measure of the union $I_n(\un t)$ of all the sets $[\un i]$
with $g(\un i)<\zeta n$ is bounded by $$ \si^{-(1-\zeta)n}
(e^\gamma p^\zeta  q)^n. $$ Since $\si>q$ we may choose
$\zeta$ so small that $e^\gamma p^\zeta q <
\si^{(1-\zeta)}$. Then $m(I_n(\un t))\le \tau^n$ with
$\tau=e^{\gamma+\zeta-1}\cdot p^\zeta \cdot q <1$ for big
enough $n\ge N$. Note that $\tau$ and $N$ do not depend on
$\un t$. Setting $$ I_n=\mcup_{\un t\in T^\NN} \big(\{\un
t\} \times I_n(\un t)\big) $$ we also have $(\th^\NN\times
m)(I_n)\le \tau^n$ for all big $n\ge N$ and for every Borel
probability $\th$ on $T$, by Fubini's Theorem. Since
$\sum_n (\th^\NN\times m)(I_n) < \infty$ then
Borel-Cantelli's Lemma implies $$ (\th^\NN\times m) \left(
\mcap_{n\ge1} \mcup_{k\ge n} I_k \right) =0 $$ and this
means that $\th^\NN\times m$ almost every $(\un t,x)\in
T^\NN\times M$ satisfies~(\ref{eq.frequency}). \cqd

\cle \label{pr.nonunifexp.localdiffeo} Let $\{
B_1,\ldots,B_p,B_{p+1},\ldots,B_{p+q}\}$ be a measurable
cover of $M$ satisfying conditions 1 and 2 of
Lemma~\ref{l.frequency}. For $0<\lambda<1$ there are
$\eta>0$ and $c_0>0$ such that, if $f_t$ also satisfies for
all $t\in T$ \begin{enumerate} \item[3.] $\| Df_t(x)^{-1}
\| \le \lambda <1$ for $x\in B_1,\ldots,B_p$; \item[4.]
$\|Df_t(x)^{-1}\| \le 1+\eta$ for $x\in
B_{p+1},\ldots,B_{p+q}$; \end{enumerate} then we have
 for $f\equiv
f_{t^*}$, where $t^*$ is some given point in $T$,
\begin{equation}\label{eq.nonunifexp.sum}
\limsup_{n\to+\infty}\frac1n\sum_{j=0}^{n-1} \log\|
Df(f^j_{\un t}(x))^{-1}\| \le -c_0 \end{equation} for
$\th^\NN\times m$ almost all $(\un t,x)\in T^\NN\times M$,
where $\th$ is any Borel probability measure on $T$.
Moreover the first hyperbolic time map $h:T^\NN\times M\to
\ZZ^{+}$ satisfies $$ (\th^\NN\times m) \{ (\un t,x)\in
T^\NN\times M: h(\un t,x) > k \} \le a_k \qand
\sum_{k=1}^\infty ka_k < \infty $$ with  $(a_k)_k$
independent of the choice of $\th$. \fle

\dem Let $\zeta>0$ be the constant provided by
Lemma~\ref{l.frequency}. We fix $\eta>0$ sufficiently small
so that $\lambda^\zeta(1+\eta)\le e^{-c_0}$ holds for some
$c_0>0$ and take $(\un t,x)$
satisfying~(\ref{eq.frequency}). Conditions 3 and 4 now
imply \begin{equation}\label{eq.nonunifexp.prod}
\prod_{j=0}^{n-1}\|Df(f^j_{\un t}(x))^{-1}\| \le
\lambda^{\zeta n} (1+\eta)^{(1-\zeta)n} \le e^{-c_0 n}.
\end{equation} for large enough $n$. This
means~(\ref{eq.nonunifexp.prod}) holds for $\th^\NN\times
m$ almost every $(\un t,x)\in T^\NN\times M$.

We observe that if $h(\un t,x)=k$, then $1\le n <  k$
cannot be hyperbolic times for $(\un t,x)$. Hence $(\un
t,x)\in I_n$ for all $n=1,\ldots,k-1$. In particular $$
(\th^\NN\times m)\{ (\un t,x)\in T^\NN\times M: h(\un
t,x)=k \} \le (\th^\NN\times m)(I_{k-1}) \equiv a_k $$ and
$\sum_k k a_k\le\sum_k k\tau^{k-1} < \infty$. \cqd

Now we will show that families of $C^2$ maps satisfying
conditions 1 through 4 of Lemmas~\ref{l.frequency} and
\ref{pr.nonunifexp.localdiffeo} contain open sets of
families in the $C^2$ topology. Let $M$ be a
$n$-dimensional torus $\TT^n$ and $f_0:M\to M$ a uniformly
expanding map: there exists $0<\lambda <1$ such that
$\|Df_0(x)v\|\ge\lambda^{-1}\|v\|$ for all $x\in M$ and
$v\in T_x M$. Let also $W$ be some small compact domain in
$M$ where $f_0\mid W$ is injective. Observe that $f_0$ is a
volume expanding local diffeomorphism due to the uniform
expansion.

Modifying $f_0$ by an isotopy inside $W$ we may obtain a
map $f_1$ which coincides with $f_0$ outside $W$, is volume
expanding in $M$, i.e., $| \det Df_1(x) | >1$ for all $x\in
M$, and has bounded contraction on $W$ near 1: $\| Df_1(x)
^{-1} \| \le 1+\eta$  for every $x\in W$ and some $\eta>0$
small. This new map $f_1$ may be taken $C^1$ close to $f_0$
and we may consider a $C^2$ map $f_2$ arbitrarily $C^1$
close to $f_1$.

Now any map $f$ in a small enough $C^2$ neighborhood of
$f_2$ admits $\si>1$ such that $|\det Df(x)|\ge\si$ for all
$x\in M$ and, for $x$ outside $W$, we have $\|
Df(x)^{-1}\|\le\lambda$. If the $C^2$ neighborhood is taken
sufficiently small then we maintain $\|Df(x)^{-1}\|\le
1+\eta$ for $x\in W$ and for some small $\eta>0$. Let us
take $B_1,\ldots,B_p, B_{p+1}=W$ a partition of $M$ into
measurable sets where the restriction $f\mid B_i$ is
injective for $i=1,\ldots,p+1$. Then any continuous family
of $C^2$ maps $\Phi:T\to C^2(M,M)$ together with a family
$(\th_\ep)_{\ep>0}$ of Borel probability measures in the
metric space $T$, satisfying $\supp(\th_\ep)\to\{t^*\}$
when $\ep\to0$ and $f_{t^*}\equiv f$, for some $t^*\in T$,
is such that $f$ is non-uniformly expanding for random
orbits and $(h_\ep)_{\ep>0}$ has uniform $L^1$-tail -- by
Lemma~\ref{pr.nonunifexp.localdiffeo} with $q=1$ and
$T=\supp(\th_\ep)$ for small enough $\ep>0$. Theorem B then
shows

\cco \label{co.nonunifexp} There are open sets ${\cal U}
\subset C^2(M,M)$
 of maps in the $C^2$ topology
such that every map $f\in {\cal U}$ is a stochastically
stable non-uniformly expanding local diffeomorphism. \fco

\subsection{Viana maps}
\label{ss.vianamaps}

In what follows we  describe the class of non-uniformly
expanding maps (with critical sets) introduced by M. Viana,
referring the reader to \cite{V}, \cite{Al} and \cite{AV}
for details. Then we show that those maps satisfy the
hypotheses of  Theorem \ref{t.finite2} and Theorem
\ref{t.sts}.

 Let $a_0\in(1,2)$ be such that
the critical point $x=0$ is pre-periodic for the quadratic
map $Q(x)=a_0-x^2$. Let $S^1=\RR/\ZZ$ and $b:S^1\rightarrow
\RR$ be a Morse function, for instance, $b(s)=\sin(2\pi
s)$. For fixed small $\alpha>0$, consider the map
 \[ \begin{array}{rccc} \hat f: & S^1\times\RR
&\longrightarrow & S^1\times \RR\\
 & (s, x) &\longmapsto & \big(\hat g(s),\hat q(s,x)\big)
\end{array}
 \]
 where $\hat g$ is the uniformly expanding
map of the circle defined by $\hat{g}(s)=ds$ (mod $\ZZ$)
for some $d\ge16$, and $\hat q(s,x)=a(s)-x^2$ with
$a(s)=a_0+\al b(s)$. It is easy to check that for $\al>0$
small enough there is an interval $I\subset (-2,2)$ for
which $\hat f(S^1\times I)$ is contained in the interior of
$S^1\times I$. Thus, any map $f$ sufficiently close to
$\hat f$ in the $C^0$ topology has $S^1\times I$ as a
forward invariant region. We consider from here on these
maps $f$ close to $\hat f$ restricted to $S^1\times I$.

Taking into account the expression of $\hat f$ it is not
difficult to check that $\hat f$ (and any map $f$  close to
$\hat f$ in the $C^2$ topology) behaves like a power of the
distance close to the critical set.

\subsubsection{Non-uniform expansion}
\label{ss.non-uniform}

 The results in
\cite{V} show that if the map $f$ is sufficiently close to
$\hat f$ in the $C^3$ topology then $f$ has two positive
Lyapunov exponents almost everywhere:
 there is a constant $\lambda>0$ for which
$$\liminf_{n\rightarrow+\infty}\frac{1}{n}\log
\|Df^n(s,x)v\|\geq \lambda$$ for Lebesgue almost every
$(s,x)\in S^1\times I$ and every non-zero  $v\in
T_{(s,x)}(S^1\times I)$. In fact, the method used for
showing this result also gives that $f$ is a non-uniformly
expanding map and its orbits have slow approximation to the
critical set, as we now explain. For the sake of clearness,
we start by assuming that  $f $ has the special form
 \begin{equation}
f(s,x)=(g(s),q(s,x)),\quad\mbox{with}\quad\partial_xq(s,x)=0
\quad\mbox{if and only if}\quad x=0, \label{assumption}
\end{equation}
 and describe how the conclusions in \cite{V}
are obtained for each $C^2$ map $f$ satisfying
 \begin{equation}
\|f-\hat f\|_{C^2}\leq\al\quad\mbox{on}\quad S^1\times
I.\label{alfa}
 \end{equation}
 Then
we explain how these conclusions extend to the general
case, using the existence of a central invariant foliation,
and we show how the results in \cite{V} give the
 non-uniform expansion and slow approximation of orbits to the critical
 set for each map $f$ as in (\ref{alfa}).

 The estimates on the derivative rely on a
statistical analysis of the returns of orbits to the
neighborhood $S^1\times (-\sqrt{\al},\sqrt{\al}\,)$ of the
critical set $\cc=\{(s,x):x=0\}$.
 We set
 $$J(0)=I\setminus(-\sqrt\al,\sqrt\al)\qand
  J(r)=\{x\in I:|x|<e^{-r}\}\quad\mbox{for
$r\ge0$}.$$
 From here on we only consider points $(s,x)\in S^1\times I$
whose orbit does not hit the critical set $\cc$. This constitues
no restriction in our results, since the  set of those
points has full Lebesgue measure.

 For each integer $j\geq 0$  we
define $(s_j,x_j)=f^j(s,x)$ and
  $$r_j(s,x)=\min\left\{r\ge0:x_j\in
J(r)\right\}. $$
 Consider, for some small constant $0<\eta<1/4$,
 $$
 G=\bigg\{0\le j<n:
 r_j(s,x)\ge\bigg(\frac12-2\eta\bigg)\log\frac1\al
 \bigg\}.
 $$
Fix some integer $n\ge1$ sufficiently large (only depending
on $\al>0$).
The results in \cite{V} show that if we take
 $$
 B_2(n)=\big\{(s,x):\:\mbox{there is $1\le j<n$ with  $x_j\in J([\sqrt n])$\:}\},
 $$
where $[\sqrt n]$ is the integer part of $\sqrt n$, then we
have
 \begin{equation}
 \label{e.b2}
 m(B_2(n))\le\mbox{const}\,e^{-\sqrt n/4}
 \end{equation}
  and, for
every small $c>0$ (only depending on the quadratic map
$Q$),
 \begin{equation}
 \label{e.principal}
 \log\prod_{j=0}^{n-1}|\partial_xq(s_j,x_j)|\ge
 2cn-\sum_{j\in G}r_j(s,x)\quad\mbox{for}\quad (s,x)\notin B_2(n),
 \end{equation}
 see \cite[pp. 75 \& 76]{V}.
Moreover, if we define for $\gamma>0$
 $$
 B_1(n)=\bigg\{(s,x)\notin B_2(n):\sum_{j\in G}r_j(s,x)\ge
 \gamma n\bigg\},
 $$
then, for small $\gamma>0$, there is a constant $\xi>0$ for which
 \begin{equation}
 \label{e.b1}
 m\big(B_1(n)\big)\le e^{-\xi n},
 \end{equation}
 see~\cite[p. 77]{V}.
Taking into account the definitions of $J(r)$ and $r_j$, this shows that if we take
$\delta=(1/2-2\eta)\log(1/\al)$, then
 \begin{equation}
 \nonumber
 \sum_{j=0}^{n-1}-\log\dist_\delta(f^j(x),\cc)\leq \gamma n
 \quad\mbox{for}\quad(s,x)\notin B_1(n)\cup B_2(n).
 \end{equation}
This in particular gives that almost all orbits have slow
approximation to $\cc$.

 On the other hand, we have  for $(s,x)\in
S^1\times I$
 \begin{equation}\label{e.forma}
  \big( Df(s,x)\big)^{-1}=
 \frac{1}{\partial_x q(s,x)\partial_s g(s)}
  \left( \begin{array}{cc}
     \partial_x q(s,x) & 0 \\
   -\partial_{s}q(s,x) &  \partial_s g(s)
 \end{array} \right).
 \end{equation}
Since all the norms are equivalent in finite dimensional
Banach spaces, it no restriction for our purposes
if we take the norm of $\big(Df(s,x)\big)^{-1}$ as the
maximum of its entries. From (\ref{assumption}) and (\ref{alfa})
we deduce that for small $\alpha$
 $$
 |\partial_s g|\geq d-\alpha,\quad |\partial_s q|\leq
 \alpha|b^\prime|+\alpha\leq 8\alpha\quad
 \mbox{and}\quad|\partial_xq|
 \leq |2x|+\alpha\leq 4,
 $$
which together with (\ref{e.forma}) gives
 $$
 \big\|\big(Df(s,x)\big)^{-1}\big\|
 =
 |\partial_xq(s,x)|^{-1},
 $$
as long as $\al>0$ is taken sufficiently small. This
implies
 \begin{equation}
 \label{e.soma}
 \sum_{j=0}^{n-1}\log\|Df(s_j,x_j))^{-1}\|=-
 \sum_{j=0}^{n-1}\log|\partial_xq(s_j,x_j)|
 \end{equation}
for every $(s,x)\in S^1\times I$.
If we choose $\gamma<c$, then we have
 \begin{equation}
 \label{e.soma2}
 \sum_{j=0}^{n-1}\log|\partial_xq(s_j,x_j)|=
 \log\prod_{j=0}^{n-1}|\partial_xq(s_j,x_j)|\ge cn
 \end{equation}
for every $(s,x)\notin B_1(n)\cup B_2(n)$ (recall (\ref{e.principal}) and
the definition of $B_1(n)$).
We conclude from (\ref{e.soma}) and (\ref{e.soma2}) that
 $$
 \sum_{j=0}^{n-1}\log\|Df(s_j,x_j))^{-1}\|\le
 -cn
 \quad
 \mbox{for}\quad(s,x)\notin B_1(n)\cup B_2(n),
 $$
which, in view of the estimates on the Lebesgue measure of
$B_1(n)$ and $B_2(n)$,  proves that $f$ is a non-uniformly
expanding map.

\medskip

 Now we describe how in \cite{V} the
same conclusions are obtained  without assuming
(\ref{assumption}).
Since $\hat f$ is strongly expanding in the horizontal
direction, it follows from the methods of \cite{HPS} that
any map $f$ sufficiently close to $\hat f$ admits a unique
invariant central foliation $\cf^c$ of $S^1\times I$ by
smooth curves uniformly close to vertical segments, see
\cite[Section 2.5]{V}. Actually, $\cf^c$ is obtained as
the set of integral curves of a vector field $(\xi^c,1)$
in $S^1\times I$
with $\xi^c$ uniformly close to zero.
 The previous analysis can then be
carried out in terms of the expansion of $f$ along this
central foliation $\cf^c$. More precisely, $|\partial_x
q(s,x)|$ is replaced by
 $$
|\partial_c q(s,x)| \equiv |Df(s,x)v_c(s,x)|,
 $$ where
$v_c(s,x)$ is a unit vector tangent to the foliation at
$(s,x)$. The previous observations imply that $v_c$ is
uniformly close to $(0,1)$ if $f$ is close to $\hat f$.
Moreover, cf. \cite[Section 2.5]{V}, it is no restriction
to suppose $|\partial_c q(s,0)|\equiv 0$, so that
$\partial_c q(s,x)\approx |x|$, as in the unperturbed case.
Indeed, if we define the {\em critical set} of $f$ by
 $$
 \cc=\{(s,x)\in S^1\times I:\partial_c q(s,x)=0\}.
 $$
 by an easy
implicit function argument it is shown in  \cite[Section
2.5]{V} that $\cc$  is the graph of some $C^2$ map $\eta :
S^1\rightarrow I$ arbitrarily $C^2$-close to zero if
$\alpha$ is small. This means that up to a change of
coordinates $C^2$-close to the identity we may suppose that
$\eta\equiv 0$ and, hence, write for $\al>0$ small
 $$
 \partial_c q(s,x)=x\psi(s,x)\quad\mbox{with $|\psi+2|$ close
 to zero}.
 $$
This provides an analog to the second part of assumption
(\ref{assumption}).
 At this point, the
arguments  apply with $\partial_xq(s,x)$ replaced by
$\partial_c q(s,x)$, to show that orbits have slow
approximation to the critical set $\cc$ and
 $
 \prod_{i=0}^{n-1} |\partial_c q(s_i,x_i)|
 $ grows
exponentially fast for Lebesgue almost every $(s,x)\in
S^1\times I$. A matrix formula for
$\big(Df^n(s,x)\big)^{-1}$ similar to that in (\ref{e.forma}) can
be obtained if we replace the vector $(0,1)$ in the
canonical basis of the  space tangent to $S^1\times I$ at $(s,x)$ by
$v_c(s,x)$, and consider the matrix of $\big(Df^n(s,x)\big)^{-1}$
with respect to the new
basis.

\medskip

For future reference, let us make some considerations on
the way the sets $B_1(n)$ and $B_2(n)$ are obtained. Let
$X:S^1\ra I$ be a smooth map whose graph in $ S^1\times I$
is nearly horizontal (see the notion of admissible curve in
\cite[Section~2]{V} for a precise definition). Denote $\wh
X_n(s)=f^n\big(s,X(s)\big)$ for $n\ge0$ and $s\in S^1$.
Take
 some leaf $L_0$ of the foliation  $\cf^c$.
Letting $L_n=f^n(L_0)$ for $n\ge 1$, we define a sequence
of Markov partitions $(\cp_n)_n$ of $S^1$ in the following
way:
  $$
   \cp_n=\left\{[s^\prime,s^{\prime\prime})\colon (s^\prime,s^{\prime\prime})
   \mbox{ is a connected component of }
   \widehat{X}_n^{-1}\big((S^1\times I)\setminus L_n\big)\right\}.
 $$
It is easy to check that $\cp_{n+1}$ refines $\cp_n$ for
each $n\ge 1$ and
 $$
 (d+\mbox{const}\,\al)^{-n}\leq |\omega|\leq (d-\mbox{const}\,\al)^{-n}
 $$
for each $\omega\in\cp_n$. Due to the  large expansion of
$f$ in the horizontal direction, we have that if $J\subset
I$ is an interval with $|J|\le\al$, then for each
$\omega\in\cp_n$
 \begin{equation}
 \label{e.cor2.3}
 m\big(\{s\in \omega\colon \wh X_j(s)\in S^1\times
 J\}\,\big)\le \mbox{const}\,\sqrt{|J|}\,m(\omega)
 \end{equation}
see \cite[Corollary~2.3]{V}. The estimate (\ref{e.b2}) on
the Lebesgue measure of $B_2(n)$ is now an easy consequence
of (\ref{e.cor2.3}). For that we only have to compute the
Lebesgue measure of $B_2(n)$  on each horizontal line of
$S^1\times I$ and integrate. The estimate (\ref{e.b2}) on
the Lebesgue measure of $B_1(n)$ is is obtained by mean of
a large deviations argument applied to the horizontal
curves in $S^1\times I$.

\cre \label{r.uniform} The choice of the
constants $c, \xi,\gamma$ and $ \delta $ only depends on
the quadratic map $Q$ and $\al>0$. In particular the decay
estimates on the Lebesgue measure of $B_1(n)$ and $B_2(n)$
depend only  on the quadratic map $Q$ and $\al>0$.
 \fre

\subsubsection{Random perturbations}

Let $f$ be
 close to $\hat f$ in the $C^3$ topology. As we have
seen before, it is no restriction if we assume that
$\cc=\{(s,x)\in S^1\times I\colon x=0\}$ is the critical set of $f$.
Fix $\{\Phi,(\theta_\ep)_\ep\}$ a random perturbation of
 $f$ for which (\ref{e.perturbation}) holds.
Our goal now is to prove that any
such $f$ satisfies the hypotheses of
Theorems~\ref{t.finite2} and \ref{t.sts} for $\ep>0$
sufficiently small, and thus conclude
that $f$ is stochastically stable. So, we want  to show
that if $\ep>0$ is  small enough then
 \begin{itemize}
 \item $f$ is non-uniformly expanding for random orbits;
 \item  random orbits have slow approximation to the critical set
 $\cc$;
 \item the family of hyperbolic time maps $(h_\ep)_\ep$ has uniform $L^1$-tail.
 \end{itemize}

We remark that in the estimates we have obtained for
$\log\|(Df(s_j,x_j))^{-1}\|$ and
$\log\dist_\delta(x_j,\cc)$ over the orbit of a given point
$(s,x)\in S^1\times I$, we can easily replace the iterates
$(s_j,x_j)$  by random iterates $(s^j_{\un t},x^j_{\un t})=
f_{\un t}^j(s,x)$. Actually, the methods used for obtaining
estimate (\ref{e.principal}) rely on a delicate
decomposition of the orbit of a given point $(s,x)$ from
time 0 until time $n$ into finite pieces according to its
returns to the neighborhood $S^1\times
(-\sqrt\al,\sqrt\al)$ of the critical set. The main tools
are \cite[Lemma~2.4]{V} and \cite[Lemma~2.5]{V} whose
proofs may easily be mimicked for random orbits. Indeed,
the important fact in the proof of the referred lemmas is
that orbits of points in the central direction stay close
to  orbits of the quadratic map $Q$ for long periods, as long
as $\al>0$
is taken sufficiently small. Hence, such results can easily
be obtained for random orbits as long as we take $\ep>0$
with $\ep\ll\al$ and perturbation vectors $\un
t\in\supp(\theta_\ep)$.

 Thus, the
procedure of \cite{V} described in Subsection
\ref{ss.non-uniform} applies to this situation, and we are
able to prove that there is $c>0$, and for $\gamma>0$ there
is $\de>0$, such that
 $$
 \sum_{j=0}^{n-1}\log\|Df(s^j_{\un t},x^j_{\un t}))^{-1}\|\le
 -cn
 \qand
 \sum_{j=0}^{n-1}-\log\dist_\delta(x^j_{\un t},\cc)\leq \gamma n
 $$
 for $(s,x)\notin B_1(n)\cup B_2(n)$, where
 $B_1(n)$ and $ B_2(n)$ are subsets $S^1\times I$ with
 $$
 m\big(B_1(n)\big)\le e^{-\xi n}\qand
 m(B_2(n))\le\mbox{const}\,e^{-\sqrt n/4}
 $$
 for some constant $\xi>0$ depending only on $\gamma$.
 This gives the
  non-uniform expansion and slow approximation to the critical set
   for random orbits. Moreover, the arguments show that we may
   take the map $N_\ep$ with
 $$ (\th_\ep^\NN\times
 m)\left(\big\{(\un t,x)\in T^\NN\times M\colon N_\ep(\un
 t,x)>n\big\}\right)\leq \mbox{const}\,e^{-\sqrt n/4},
 $$
 thus giving that the family of first hyperbolic time maps
 has uniform $L^1$-tail (see the considerations at the beginning
 of Section \ref{s.applications}).

\medskip

For the sake of completeness, an explanation is required on
the way the Markov partitions $\cp_n$ of $S^1$ can be
defined in this case, in order to obtain the estimates on
the Lebesgue measure of $B_1(n)$ and $B_2(n)$. We consider
$M=S^1\times I$ and define the skew-product map
 \[
\begin{array}{rccc} F: & T^\NN\times M &\longrightarrow &
T^\NN\times M,\\
 & (\un t, z) &\longmapsto & \big(\sigma(\un t),f_{t_1}(z)\big)
\end{array}
 \]
 where $\si$ is the left shift map. Writing
$f_t(z)=\big(g_t(z),q_t(z)\big)$ for $z=(s,x)\in S^1\times
I$, we have that $q_t(s,\cdot)$ is a unimodal map close to
$\hat q$ for all $s\in S^1$ and $t\in \supp(\theta_\ep)$
with $\ep>0$ small.

\cpr\label{p.foliation}
 Given $\un t\in T^\NN$ there is a
 $C^1$ foliation $\cf^c_{\un t}$ of $M$ such that if $L_{\un t}(z)$
is the leaf of $\cf^c_{\un t}$ through a point $z\in M$,
then
 \begin{enumerate}
 \item  $L_{\un t}(z)$ is a $C^1$ submanifold of
 $M$ close to a vertical line in the $C^1$ topology;
 \item $f_{t_1}\big(L_{\un t}(z)\big)$ is contained in $L_{\si\un t}\big(f_{t_1}(z)\big)$.
 \end{enumerate}
 \fpr
 \dem
This will be obtained as  a  consequence of the fact that
the set of vertical lines constitutes a normally expanding
invariant foliation for $\hat f$. Let $\ch$ be the space of
continuous maps $\xi:T^\NN\times M\rightarrow [-1,1]$
endowed with the sup norm, and define the map
$A:\ch\rightarrow\ch$ by
 $$
A\xi(\un t,z)=\frac{\partial_xq_{t_1}(z)\xi( F(\un t,
z))-\partial_xg_{t_1}(z)} {-\partial_s q_{t_1}(z)\xi(F(\un
t,z) )+\partial_s g_{ t_1}(z)}, \quad \un
t=(t_1,t_2,\dots)\in T^\NN \qand z\in M.
 $$
 Note that $A$ is well-defined,
since
 $$
 |A\xi(\un t,z)|\leq \frac{(4+\al+\epsilon)+\al+\epsilon}
 {-(\mbox{const }\alpha+\epsilon )+(d-\al-\epsilon)}
 <1\label{campo}
 $$
 for small $\alpha>0$ and $\epsilon>0$. Moreover,
$A$ is a contraction on $\ch$: given $\xi,\zeta\in\ch$ and
$(\un t,z)\in T^\NN\times M$ then
 \begin{eqnarray*}
\lefteqn{|A\xi(\un t,z)-A\zeta(\un t,z)|}\hspace{2cm}\\
 &\leq &
 \frac{|\mbox{det}Df_{t_1}(z)|\cdot|\xi(\un t, z)-\zeta(\un t,z)|}
 {\big|\big(\!-\partial_s q_{t_1}(z)\xi(F(\un t,z))+\partial_s
g_{t_1}(z)\big)\cdot\big (\!-\partial_s
q_{t_1}(z)\zeta(F(\un t,z))+\partial_s
g_{t_1}(z)\big)\big|}\\ &\leq &
\frac{\big((d+\alpha+\epsilon)(4+\alpha+\epsilon)+\alpha+
\epsilon\big)\cdot\big|\xi(\un t, z)-\zeta(\un t,z)\big|}
{(d-\mbox{const}\alpha-\epsilon)^2}.
 \end{eqnarray*}
  This last
quantity can be made smaller than $|\xi(\un t,z)-\eta(\un
t, z)|/2$, as long as $\alpha$ and $\epsilon$ are chosen
sufficiently small. This shows that $A$ is a contraction on
the Banach space $\ch$, and so it has a unique fixed point
$\xi^c\in\ch$.

 It is no restriction for our purposes if we
think of $T$ as being equal to $\supp(\theta_{\ep})$ for
some small $\ep$.
Note that the map $A$ depends
continuously on $F$ and for $\ep>0$ small enough  the fixed
point of $A$ is close to the zero constant map. This holds
because we are choosing $\supp(\theta_{\ep})$ close to
$\{t^*\}$, $f_{t^*}=f$ and $f$ close to $\hat f$. Then, for
$\ep>0$ small enough, we have $\xi^c(\un t, \cdot)$ uniformly
close to $\xi^c(\un t^*, \cdot)$ and it is not hard to
check that $\xi^c_0=\xi^c(\un t^*, \cdot)$ is precisely the map
whose integral leaves of the vector field $(\xi_0^c,1)$
give the invariant foliation $\cf^c$ associated to $f_{t^*}=f$.
Since this foliation depends continuously on the dynamics
and for $f=\hat f$ we have $\xi_0^c\equiv 0$
(see \cite[Section~2.5]{V}), we finally deduce
that $\xi^c(\un t, \cdot)$ is uniformly close to zero for
small $\ep>0$.

 We have defined $A$ in such a way
that if we take $E^c(\un t, z)=\mbox{span}\{(\xi^c(\un
t,z),1)\}$, then for every $\un t\in T^\NN$ and $z\in
S^1\times I$
 \begin{equation} Df_{t_1}(z)E^c(\un t,z)\subset
E^c(F(\un t,z)).\label{ec} \end{equation}
 Now, for fixed $\un t\in T^\NN$, we take $\cf^c_{\un t}$ to
be the set of integral curves of the vector field
$z\rightarrow (\xi^c(\un t ,z),1)$ defined on $S^1\times
I$. Since the vector field is taken of class $C^0$, it does
not follow immediately that through each point in
$S^1\times I$ passes only one integral curve. We will prove
uniqueness of solutions by using the fact that the map $f$
has a big expansion in the horizontal direction.

 Assume, by
contradiction, that there are two distinct integral curves
$Y, Z\in\cf^c_{\un t}$ with a common point. So we may take
three distinct nearby points $z_0,z_1,z_2\in S^1\times I$
such that $z_0\in Y\cap Z$, $z_1\in Y$, $z_2\in Z$ and
$z_1,z_2$ have the same $x$-coordinate. Let $X$ be the
horizontal curve joining $z_1$ to $z_2$. If we consider
$X_n=\pi_2\circ F^n(\un t,X)$ for $n\ge 1$, where $\pi_2$
is the projection from $T^\NN\times S^1\times I$ onto
$S^1\times I$, we have that the curves $X_n$ are nearly
horizontal and grow in the horizontal direction (when $n$
increases) by a factor close to $d$ for small $\al$ and
$\ep$, see \cite[Section 2.1]{V}. Hence, for large $n$,
$X_n$ wraps many times around the cylinder $S^1\times I$.
On the other hand, since $Y_n=\pi_2\circ F^n(\un t,Y)$ and
$Z_n=\pi_2\circ F^n(\un t,Z)$ are always tangent to the
vector field $z\rightarrow \big(\xi^c(\sigma^n\un
t,z),1\big)$ on $S^1\times I$, it follows that all the
iterates of $Y_n$ and $Z_n$ have small amplitude in the
$s$-direction. This gives a contradiction, since the closed
curve made by $Y$, $Z$ and $X$ is homotopic to zero in
$S^1\times I$ and the closed curve made by $Y_n$, $Z_n$ and
$X_n$ cannot be homotopic to zero for large $n$. Thus, for
fixed $\un t\in T^\NN$ we have uniqueness of solutions of
the vector field $z\rightarrow (\xi^c(\un t,z),1)$, and
from (\ref{ec}) it follows that $\cf^c_{\un t}$ is an
$F$-invariant foliation of $M$ by nearly vertical leaves.
\cqd

Now, using the foliations given by the previous proposition
we are also able to define the Markov partitions of $S^1$
in this setting.  Given any smooth map $X:S^1\ra I$ whose
graph is nearly horizontal,  denote
 $\wh X_{\un t}^n
(s)=f_{\un t}^n\big(s,X(s)\big)$ for $n\ge0$ and $s\in
S^1$. Take
 some leaf
  $L_{\un t}^0$ of the foliation $\cf^c_{\un t}$.
Letting $L^n_{\un t}=f^n_{\un t}(L_{\un t})$ for $n\ge 1$,
we define the sequence of Markov partitions $(\cp^n_{\un t
})_n$ of $S^1$ as
  $$
  \cp^n_{\un t}=\left\{[s^\prime,s^{\prime\prime})\colon
 (s^\prime,s^{\prime\prime})
   \mbox{ is a connected component of }
 (\widehat{X}^n_{\un t})^{-1}\big((S^1\times I)\setminus L^n_{\un t
}\big)\right\}.
 $$
It is easy to check that $\cp^{n+1}_{\un t }$ refines
$\cp^n_{\un t }$ for each $n\ge 1$ and, taking $\ep\ll\al$,
 $$
 (d+\mbox{const}\,\al)^{-n}\leq |\omega|\leq (d-\mbox{const}\,\al)^{-n}
 $$
for each $\omega\in\cp^n_{\un t }$. This permits to obtain
estimates for the Lebesgue measure of the sets $B_1(n)$ and
$B_2(n)$ exactly in the same way as before also with the
constants only depending on the quadratic map $Q$ (cf.
Remark \ref{r.uniform}).


\bigskip

\noindent Jos\'e Ferreira Alves ({\tt jfalves@fc.up.pt}) \\
V\'\i tor Ara\'ujo ({\tt vdaraujo@fc.up.pt})\\
Centro de Matem\'atica da Universidade do Porto\\
Pra\c ca Gomes Teixeira,
4099-002 Porto, Portugal\\

\end{document}